\documentclass[10pt,conference]{IEEEtran}

\usepackage[cmex10]{amsmath}
\usepackage{amsfonts,amssymb,graphicx}
\graphicspath{ {./} }
\DeclareGraphicsExtensions{.pdf}
\usepackage{algorithm,algpseudocode}
\usepackage{caption}
\usepackage{subcaption}
\usepackage[mathscr]{eucal}

\usepackage{mathtools}
\usepackage[subnum]{cases}

\usepackage{float}

\usepackage{url}

\usepackage{verbatim}

\usepackage{booktabs}
\usepackage{multirow} 

\usepackage{fixltx2e}

\newcommand{\Abs}[1]{\ensuremath{\left\vert #1 \right\vert}}
\newcommand{\InfNorm}[1]{\ensuremath{ \Norm{#1}_{\infty}}}
\newcommand{\TwoNorm}[1]{\ensuremath{ \Norm{#1}_{2}}}

\newcommand{\FNorm}[1]{\ensuremath{ \Norm{#1}_{F}}}
\newcommand{\Norm}[1]{\ensuremath{ \left\Vert #1 \right\Vert}}

\newcommand{\Mat}[1]{\ensuremath{ \boldsymbol{\mathbf{#1}}}}
\newcommand{\Vect}[1]{\ensuremath{ \boldsymbol{\mathbf{#1}}}}

\newcommand{\R}{\ensuremath{\mathbb{R}}}
\newcommand{\Rn}{\ensuremath{\R^n}}
\newcommand{\argminy}{\ensuremath{\displaystyle\operatornamewithlimits{arg\,min}_y}}
\newcommand{\thinColon}{%
  \nobreak
  \mskip1mu 
  \mathpunct{}%
  \nonscript
  \mkern-\thinmuskip
  {:}%
  \mskip1mu
  \relax
}
\newcommand{\IndexRange}[2]{\ensuremath{{#1}\thinColon{}{#2}}}

\newcommand{\bvector}[3]{
    \left\{\!\!\!
        \begin{array}{c}
            #1     \\
            #2     \\
            #3     \\
        \end{array}
    \!\!\!\right\}
}

\title{Exploiting Data Representation for Fault Tolerance}

 \author{
 \IEEEauthorblockN{
 James Elliott\IEEEauthorrefmark{1}\IEEEauthorrefmark{2},
 Mark Hoemmen\IEEEauthorrefmark{2}, and
 Frank Mueller\IEEEauthorrefmark{1}}
 \IEEEauthorblockA{
 \IEEEauthorrefmark{1}
 Computer Science Department, North Carolina State University,
 Raleigh, NC
 }
 \IEEEauthorblockA{
 \IEEEauthorrefmark{2}
 Sandia National Laboratories
 Albuquerque, NM
 }
 }

\begin{document}

\maketitle

\begin{abstract}
We explore the link between data representation and
soft errors in dot products. We present an analytic model for the absolute
error introduced should a soft error corrupt a bit in an IEEE-754
floating-point number. We show how this finding relates to
the fundamental linear algebra concepts of normalization and matrix
equilibration. We present a case study illustrating that the
probability of experiencing a large error in a dot product is
minimized when both vectors are normalized. Furthermore, when data is
normalized we show that the absolute error is less than one or very
large, which allows us to detect large errors. We demonstrate how this
finding can be used by instrumenting the GMRES iterative solver. We
count all possible errors that can be introduced through faults in
arithmetic in the computationally intensive orthogonalization phase,
and show that when scaling is used the absolute error can be bounded
above by one.
\end{abstract}

\section{Introduction} \label{S:intro}
In the field of high-end computing (HEC) the notion of reliability has
tended to focus on keeping thousands of physical nodes operating
cooperatively for extended periods of time. As chip manufacturing and
power requirements continue to advance, soft errors are becoming more
apparent \cite{jje:Michalak:2012}. This implies that reliability
research must address the case that the machine does not crash, but
that outputs during computation may be silently incorrect. There have
been many studies into hardening numerical kernels against soft
errors, that is the researchers attempt to preserve the illusion of a
reliable machine by detecting and correcting all soft errors. We take
a more analytical approach.
Instead of focusing on detection/correction, we seek to study how the
data operated on impacts the errors that we can observe given soft
errors in data --- called silent data corruption (SDC).


The driving motivation behind our work is the uncertainty surrounding
the reliability of an exascale-class machine
\cite{elnozahy08,jje:kogge2008exascale,cappello09}.  We attempt to
avoid speculation over what hardware may be used in future (or
present) HEC deployments, and instead analyze how a single soft error
in an IEEE-754 floating-point number behaves.  It has already been
shown that existing and decommissioned HEC deployments have suffered
from SDC \cite{jje:Michalak:2012,jje:haque2010hard}. For the prior
reasons, we seek to study the link between the data operated on and
soft errors.  We intentionally perform our research subject to the
IEEE 754 specification, which we believe will be used regardless of
the architecture. We also restrict our analysis to single bit flips.
This gives us a base line from which to draw higher-level conclusions
related to multiple bit flips, and lets us isolate the impact of a bit
flip.

IEEE 754 both defines the binary \emph{representation} of data, and
bounds the \emph{rounding error} committed by arithmetic operations.
This work focuses on data representation.  The effects of rounding
error on numerical algorithms, including those studied in this paper,
have been extensively studied; see e.g., \cite{paige2006modifed}.
However, these results generally only apply to \emph{small} errors,
such as those resulting from rounding.  Bit flips can be huge and thus
require different methods of analysis, like those presented in this
paper.
\\

\textbf{We present the following contributions}:
\begin{itemize} 
\item We model single bit upsets in IEEE-754 scalars analytically,
  and extend this modeling to dot products.
\item We demonstrate both experimentally (via Monte Carlo sampling) and
analytically that dot products performed on normalized numbers have a
significantly lower probability of experiencing large error than dot
products with values of varying magnitudes.
\item We relate our finding that normalized vectors minimize
absolute error to matrix equilibration, and correlate this
finding to two highly used numerical kernels (Gram-Schmidt
orthogonalization and the Arnoldi process).
\item We demonstrate the utility of our finding by instrumenting the
Generalized Minimum Residual Method (GMRES). We show that for the
dot product intensive orthogonalization kernel, we can restrict errors
arising from single bit upsets to being less than one, or being very
large and easily detected.
\item We articulate how studying single bit flips enables us to draw
conclusions about multiple bit upsets.
\end{itemize}

\section{Related Work} \label{S:related_work}

Researchers have approached the problem of SDC in numerical algorithms
in various ways.  Many take the approach of treating an algorithm as a
black box and observing the behavior of these codes when run with soft
errors injected. Recently,
\cite{jje:iterative:Howle:2012,jje:iterative:Howle:2010} analyzed the
behavior of various Krylov methods and observed the variance in
iteration count based on the data structure that experiences the bit
flip. Shantharam et al.\ \cite{jje:iterative:Shantharam:2011} analyzed
how bit flips in a sparse matrix-vector multiply (SpMV) impact the
$L^2$ norm and observe the error as CG is run.  Bronevetsky et al.\
\cite{jje:iterative:Bronevetsky:2008, jje:iterative:Sloan:2012}
analyzed several iterative methods documenting the impact of randomly
injected bit flips into specific data structures in the algorithms and
evaluated several detection/correction schemes in terms of overhead
and accuracy.  Exemplifying the concept of black-box analysis,
\cite{jje:Li:2012} presents BIFIT for characterizing applications
based on their vulnerability to bit flips.  Rather than focusing on
how to preserve the illusion of a reliable machine or devising a
scheme to inject soft errors, we investigate an avenue mostly ignored,
which is how the data in the algorithm can be used to mitigate the
impact of a bit flip.


Hoemmen and Heroux proposed a radically different approach. Rather
than attempt to detect and correct soft errors, they use a ``selective
reliability'' programming model to make the algorithm converge
\emph{through} soft errors \cite{jje:bridges2012fault}. Sao and Vuduc
showed that reliably restarting iterative solvers enables convergence
in the presence of soft errors \cite{sao13}.  In the same vein,
Elliott et al.\ showed that bounding the error introduced in the
orthogonalization phase of GMRES lets FT-GMRES converge with minimal
impact on time to solution \cite{elliott13gmres}.  Boley et al.\ apply
backward error analysis to linear systems, in order to distinguish
small error due to rounding from inacceptably large error due to
transient hardware faults \cite{jje:Boley:1995}.  In general, our
work complements this line of research. While Hoemmen and Sao have
investigated algorithms that can converge through error, we show that
in certain numerical kernels, the data itself can have a ``bounding''
effect. For example, coupled with \cite{elliott13gmres}, we improve
the likelihood that errors fall within the derived bound.

Algorithm-based fault tolerance (ABFT) provides an approach to detect
(and optionally correct) faults, which comes at the cost of increased
memory consumption and reduced performance
\cite{jje:abft:Huang:1984,jje:abft:Du:2012}. The ABFT work by Huang et
al. \cite{jje:abft:Huang:1984} was proven by Anfinson et
al. \cite{jje:abft:Anfinson:1988} to work for several matrix
operations, and the checksum relationship in the input checksum
matrices is preserved at the end of the computation. Consequently, by
verifying this checksum relationship in the final computation results,
errors can be detected at the end of the computation.  Recent work has
looked at extending ABFT to additional matrix factorization algorithms
\cite{jje:abft:Du:2012} and as an alternative to traditional
checkpoint/restart techniques for tolerating fail-stop failures
\cite{jje:abft:Davies:2011,jje:abft:Chen:2008:1,jje:abft:Chen:2011}.

Costs in terms of extra memory and computation required for ABFT may
be amortized for dense linear algebra, and such overheads have been
analyzed by many (e.g.,
\cite{jje:abft:Yamani:2001,jje:abft:Banerjee:1990,jje:abft:Kim:1996}).
Algorithms must be manually redesigned for ABFT support by accounting
for numerical properties (e.g., invariants).  A more fundamental
problem is that traditional checksums and error-correcting codes do
not suit floating-point computations well \cite{jje:Boley:1995}.  Such
computations naturally commit rounding error, which exact (bitwise)
codes forbid.  Inexact codes (that use floating-point sums) can be
sensitive to rounding error, and commit it themselves.  It is possible
that more expensive recovery and significantly more redundant storage
could help \cite{candes2006highly}.  However, works like
\cite{jje:bridges2012fault,sao13,elliott13gmres,jje:Boley:1995}
suggest that \emph{correcting} faults might not be necessary, if their
effects on the algorithm are detectable and bounded.  In general, this
paper favors ``opening up the black box'' and understanding the
effects of soft error on algorithms, rather than trying to detect and
correct all such errors before they affect algorithms' behavior.

\section{Project Overview} \label{S:project}

To explore the relation between data representation and soft errors,
we first construct an analytic model of a soft error in an IEEE-754
floating-point scalar, and then extend this to a dot product.  We
uncover through analysis that the binary pattern of the exponent can
be exploited for fault tolerance. We show this graphically via a case
study using Monte Carlo sampling of random vectors, and then extend
the idea of data scaling to matrices by using sparse matrix
equilibration.  To demonstrate the feasibly and utility of our work we
analyze the GMRES algorithm and instrument the computationally
intensive orthogonalization phase. We count the possible absolute
errors that can be introduced via a bit flip in a dot product, and
show that scaling data lowers the likelihood of observing large,
undetectable errors.
\\
\\
\emph{This paper is organized as follows}:
\begin{enumerate}
\item In Section~\ref{S:fault_model}, we construct an analytic model
of the absolute error for single bit upsets in IEEE-754 floating-point
numbers.
\item In Section~\ref{S:fault_model_eval}, we extend our model of
faults in IEEE-754 scalars to vectors of arbitrary values, and present
examples of how data scaling impacts the binary representation and
absolute error we can observe.
\item In Section~\ref{S:case_study}, we perform a case study using
Monte Carlo sampling of 10,404,000,000 random vectors of various
magnitudes, and graphically show how the error is minimized when
operating on values less than 1.
\item In Section~\ref{S:gmres}, we link data scaling to sparse matrix
equilibration, and instrument and evaluate the impact of a soft error
in the computationally intensive orthogonalization phase of GMRES.
\end{enumerate}

\section{Fault Model} \label{S:fault_model}
The premise of our work is that a silent, transient bit flip impacts
data.
Before we can perform any analysis or experimental work, we must
define how such a bit flip would impact an algorithm, and how we
enforce that the bit flip was transient. To achieve this goal, we
build our model around the basic concept that when an algorithm uses data,
this translates into some set of operations being performed
on the data. Should a bit flip perturb our data, some operation will
use a corrupt value, rather than the correct value.
The output of this single operation will then contain a tainted value,
and this tainted value could cause the solution to be incorrect. Note
that a transient bit flip may cause a persistent error in the output
depending on how the value is used.



A side benefit of an operation-centric model is that we naturally avoid
a pitfall to which arbitrary memory fault injection succumbs, namely
that if a bit flip impacts data (or memory) that is never used (read)
then this fault \emph{cannot lead to a failure}. Our fault model
allows a bit flip to perturb the input to an operation performed on
the data, while not persistently tainting the storage of the inputs.
This mimics how a transient bit flip would manifest itself, e.g.,
during ALU activities.  As a result, the data that experiences the bit
flip need not show signs that it was perturbed. This model allows us
to observe the impact of transient flips on the inputs, which results
in sticky or persistent error in the result. We then utilize
mathematical analysis to model how this persistent error propagates
through the algorithm.

\subsection{Fault Characterization via Semantic Analysis}
\label{jje:sec:fault_model:characterization}
To derive a fault
model we must first understand what a fault is. Since floating-point
numbers approximate real numbers and most numerical algorithms use
real numbers, we start from the definition of a real-valued scalar
$\gamma \in \R$. The range of \emph{possible} values that $\gamma$ can
take is
\begin{equation*}
\gamma \in [-\infty,+\infty].
\end{equation*}
We assume that the IEEE-754 specification for double-precision numbers, called
\emph{Binary64}, is used to represent these numbers.  This means that $\gamma$ can 
take a fixed set of numeric values, and these values lie in
the range
\begin{equation*}\label{jje:fault_model:char:base10_noabs}
\gamma \in [ -1.80 \times 10^{308}, + 1.80 \times 10^{308}],
\end{equation*}
or using base two for the exponent
\begin{equation*}\label{jje:fault_model:char:base2_noabs}
\gamma \in [- 1.\bar{9} \times 2^{1023}, + 1.\bar{9} \times 2^{1023}],
\end{equation*}
where $1.\bar{9}$ indicates the largest possible fractional component,
and $1.0$ indicates the smallest fractional component.
A more informative range is that of $\Abs{\gamma}$, excluding $0$ and
denormalized numbers,
\begin{equation}\label{jje:fault_model:char:base10}
\Abs{\gamma} \in [ 2.23 \times 10^{-308}, 1.80 \times 10^{308}],
\end{equation}
and in semi base two
\begin{equation}\label{jje:fault_model:char:base2}
\Abs{\gamma} \in [1.0 \times 2^{-1022}, 1.\bar{9} \times 2^{1023}].
\end{equation}

\begin{figure}
    \centering
    \includegraphics[width=0.4\textwidth]{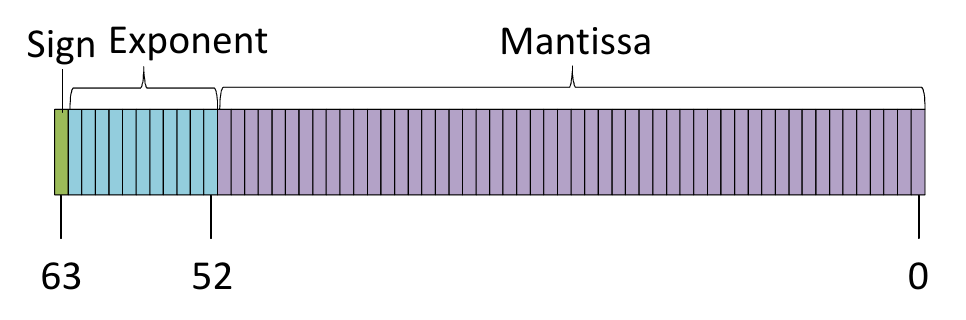}
    \caption{Graphical representation of data layout in the IEEE-754 Binary64
    specification.}
    \label{F:Binary64}
\end{figure}

To approximate real
numbers, \emph{Binary64} uses 64 bits, of which 11 are devoted to the exponent,
52 for the fractional component (we refer to as the mantissa), and one bit for
the sign. Figure~\ref{F:Binary64} shows how these bits are laid out.
In addition to numeric values, Binary64 includes two non-numeric values,
Not-a-Number (NaN) and Infinity (Inf), which may be signed to account
for infinity and values that result in undefined operations, e.g., division by zero.
The range of values in
Equations~\eqref{jje:fault_model:char:base10}~and~\eqref{jje:fault_model:char:base2}
is not continuous and has non-uniform gaps due to the discrete precision, which
is a consequence of having a fixed number of bits in the fractional component.

We can further discretize the range of possible values by recognizing
that there is a finite number of exponents that are
possible given IEEE-754 double precision, e.g.,
\begin{align*}
\gamma \in \{&	0, \pm \mathrm{Inf}, \pm \mathrm{NaN},
			\pm 2^{-1022} \times 1.x,
			\pm 2^{-1021} \times 1.x,\\
			&\dotsc, \pm 2^{0} \times 1.x, \dotsc,
			\pm	2^{1023} \times 1.x\},
\end{align*}
where $1.x$ indicates some fractional component.

Analytically, this is expressed as
\begin{equation}
\gamma = (-1)^{sign}
	\left(1 + \sum\limits_{i=0}^{51} b_i 2^{i-52}\right)
	\times 2^{e-1023},
\label{jje:eq:model:binary64}
\end{equation}
for IEEE-754 \emph{Binary64}. Note, the specification does not include a sign
bit for the exponent. Rather, IEEE floating point numbers utilize a \emph{bias}
to allow the exponent to be stored without a sign bit, which we will later
exploit for fault-resilience. Another important characteristic that
stems from the general approach of expressing numbers in exponential notation
is that we can characterize numbers by their order of magnitude. Of particular
interest is the following relation:
\begin{align}
|2^{-1022}| &\leq | 2^{-1022} \times 1.x | \nonumber \\
< |2^{-1021}| &\leq  |2^{-1021} \times 1.x| < \dots \nonumber \\
< |2^{0}| &\leq |2^{0} \times 1.x| < \dots \nonumber \\
< |2^{1023}| &\leq  |2^{1023} \times 1.x|.
\label{jje:eq:magnitude_characterizes_mantissa}
\end{align}
This means that we can use the next order of magnitude as an upper bound for
errors in the fractional component of a number --- which is practically achieved
by incrementing the exponent or multiplying by two.
We can also analytically model the number of fractional bits that could
contribute an error larger than some tolerance, since the error that \emph{could}
arise from each mantissa bit is relative to the exponent of the number. This
final step is necessary since the fractional term can take values in the range
$[1,2)$, where the left parenthesis indicates that $2$ is not a member of this
interval. We can also characterize the error that a perturbed sign bit can
contribute, and, like the fractional component, this error is relative to the
exponent of the number. Suppose the sign is perturbed in a scalar $\gamma$,
then we have $\tilde{\gamma} = -\gamma$, the absolute error is
$|\gamma - \tilde{\gamma}| = |\gamma - (-\gamma)| = 2\gamma$. This means we can
bound the error from a sign bit perturbation by incrementing the exponent of the
resulting value.

In summary, we have demonstrated that errors in IEEE-754 floating point numbers
can be characterized using the exponent of the numbers. This property allows
us to reduce the number of bits we need to consider in a fault model, since we know
that a large number of errors are bounded by the relatively small set of
possible exponents.

\subsection{Fault Characteristics of Perturbed Exponents}
In the context of IEEE-754 double precision numbers and silent data corruption,
we do not model the exponents directly. Instead, we model the biased exponents,
as they are the interesting portion of the \emph{data} that allows us to characterize
the errors that the majority of the bits present in the data can produce. For
instance, in double precision data we can characterize the errors from 53 of the
64 bits using our approach. This type of fault-characterization is impossible if
bit flips are injected randomly into the data's memory, as that approach loses
the semantic information that is implicitly present in the data.

The Binary64 specification does not store exponents directly,
instead it uses a bias of $1023$. From
\S~\ref{jje:sec:fault_model:characterization} this means we can
characterize \emph{all} faults in double precision data by analyzing
perturbations to the possible biased exponents
\begin{equation*}
	\{0,1,2,\dotsc, 1023, \dotsc, 2046\}.
\end{equation*}
Note that zero is not a biased exponent and has special meaning. In IEEE-754, a
zero pattern in the exponent with zeros in the mantissa is used to represent the
scalar zero, while a non-zero pattern in the mantissa is used to represent subnormal numbers.
We also assume the user does not perform computation on the two non-numeric
values NaN and Inf, which are represented using the biased exponent $2047$ (all
ones).
We do include zero in our analysis because it is a valid real number.

Since we are concerned with bit perturbations in the exponent, we
express the biased exponents in their binary form, e.g., 11-bit unsigned
integers presented in binary.
\begin{figure}[htp]\centering
\begin{equation*}
\begin{array}{ccccc}
\bvector{2^{-1}}{2^0}{2^{1}} & \Rightarrow & 
\bvector{1022}{1023}{1024}   & \Rightarrow &
\bvector{01111111110}{01111111111}{10000000000} \\
\mathrm{Exponent} & \qquad & \mathrm{Biased} & \qquad &
\mathrm{Storage}
\end{array}
\end{equation*}
\caption[labelInTOC]{Relation of exponent, IEEE-754 double precision bias, and
what data are actually stored.}
\label{jje:fig:exp_bias_binary}
\end{figure}
We can further expand Figure~\ref{jje:fig:exp_bias_binary} to show the potential
change to the original exponent should a bit flip occur, which will form the
basis for our fault model and analytic models.

\begin{figure*}[htbp]
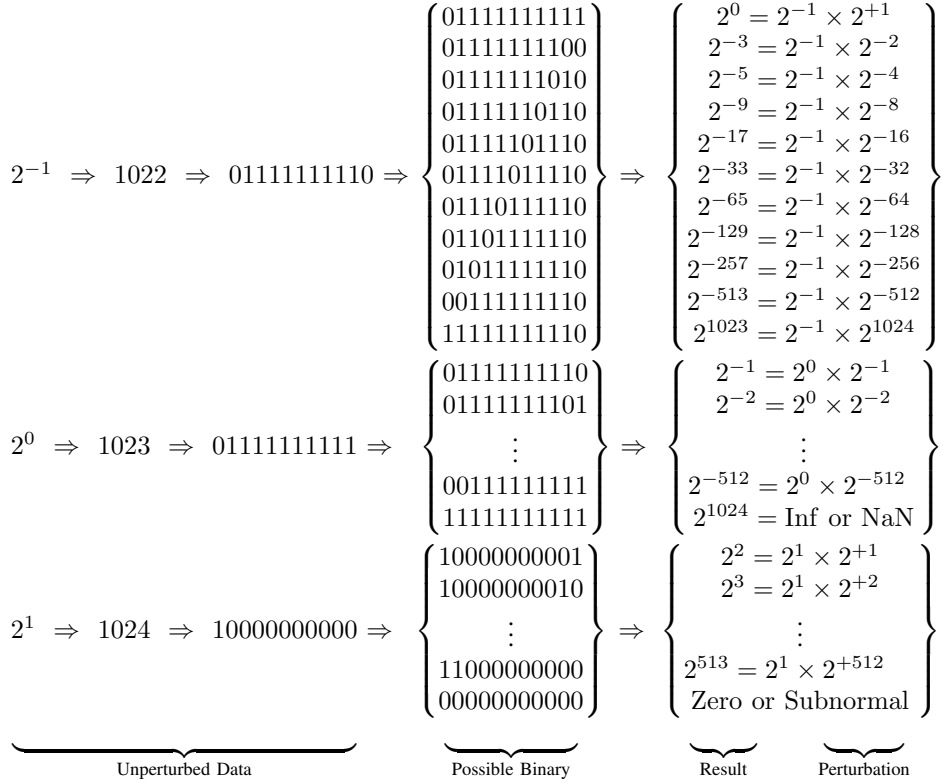

\centering
\begin{align*}
&\begin{array}{ccccccccc}
2^{-1} & \!\!\Rightarrow\!\! & 1022 & \!\!\Rightarrow\!\! & 01111111110
	& \!\!\!\!\Rightarrow\!\!\!\! &
	\left\{\!\!\!\begin{array}{c}
		01111111111 \\
		01111111100 \\
		01111111010 \\
		01111110110 \\
		01111101110 \\
		01111011110 \\
		01110111110 \\
		01101111110 \\
		01011111110 \\
		00111111110 \\
		11111111110
	\end{array}\!\!\!\right\}
	& \!\!\!\!\Rightarrow\!\! &
	\left\{\!\!\!\begin{array}{c}
		2^{0}  = 2^{-1} \times 2^{+1} \\
		2^{-3} = 2^{-1} \times 2^{-2} \\
		2^{-5} = 2^{-1} \times 2^{-4} \\
		2^{-9} = 2^{-1} \times 2^{-8} \\
		2^{-17} = 2^{-1} \times 2^{-16} \\
		2^{-33} = 2^{-1} \times 2^{-32} \\
		2^{-65} = 2^{-1} \times 2^{-64} \\
		2^{-129} = 2^{-1} \times 2^{-128} \\
		2^{-257} = 2^{-1} \times 2^{-256} \\
		2^{-513} = 2^{-1} \times 2^{-512} \\
		2^{1023} = 2^{-1} \times 2^{1024} \\
	\end{array}\!\!\!\right\}
\end{array}
\\ 
&\begin{array}{ccccccccc}
2^{0} & \!\!\Rightarrow\!\! & 1023 & \!\!\Rightarrow\!\! & 01111111111
	& \!\!\!\!\Rightarrow &
	\left\{\!\!\!\begin{array}{c}
		01111111110 \\ 
		01111111101 \\
		\vdots		\\
		00111111111 \\
		11111111111
	\end{array}\!\!\!\right\}
	& \!\!\!\!\Rightarrow\!\! &
	\left\{\!\!\!\begin{array}{c}
		2^{-1}  = 2^{0} \times 2^{-1} \\ 
		2^{-2} = 2^{0} \times 2^{-2} \\
		\vdots		\\
		2^{-512} = 2^{0} \times 2^{-512\phantom{n}} \\
		2^{1024} = \mathrm{Inf~or~NaN}
	\end{array}\!\!\!\right\}
\end{array}
\\
&\begin{array}{ccccccccc}
2^{1} & \!\!\Rightarrow\!\! & 1024 & \!\!\Rightarrow\!\! & 10000000000
	& \!\!\!\!\Rightarrow\! &
	\left\{\!\!\!\begin{array}{c}
		10000000001 \\ 
		10000000010 \\
		\vdots		\\
		11000000000 \\
		00000000000
	\end{array}\!\!\!\right\}
	& \!\!\!\!\Rightarrow\!\!\! &
	\left\{\!\!\!\begin{array}{c}
		2^{2}  = 2^{1} \times 2^{+1} \\ 
		2^{3}  = 2^{1} \times 2^{+2} \\
		\vdots		\\
		2^{513} = 2^{1} \times 2^{+512\phantom{nm}} \\
		\mathrm{Zero~or~Subnormal}
	\end{array}\!\!\!\right\} \\
\multicolumn{5}{c}{\underbrace{\hspace*{130pt}}_{\text{Unperturbed Data}}} & &
\underbrace{\hspace*{50pt}}_{\text{Possible Binary}} & &
\underbrace{\hspace*{25pt}}_{\text{Result}} \phantom{=2^{1}}
\underbrace{\hspace*{30pt}}_{\text{Perturbation}}
\end{array}
\end{align*}
\caption{Examples of how a bit flip can impact an exponent
represented using the IEEE-754 Binary64 specification.}
\label{jje:fig:exp_bias_binary_perts}
\end{figure*}

In the context of bit flips, we can view a bit flip as adding or
subtracting from the biased exponent, which in turn translates to
multiplying or dividing the number by some power of two. We model the
impact of a bit flip in the exponent as the original scalar being
magnified or minimized by a specific power of two. We illustrate this
in Figure~\ref{jje:fig:exp_bias_binary_perts}, where reading
left-to-right, we have some initial exponent, which is represented
using a bias of $1023$. The biased exponent translates to a discrete
binary pattern.
We consider all single bit flips in this binary pattern and compute
the actual perturbed exponent. Note, that the perturbation can be
modeled independent of the original exponent.

By characterizing the error introduced, we recognize that
\emph{all} mantissa bit flips introduce error that has the same
exponent as the original number, and a sign bit flip introduces error
that is only one order of magnitude larger than the original number.
Furthermore, the exponent bits can either introduce large error, or a
bit flip introduces error roughly equivalent to the order of magnitude
of the original number.

Suppose we can enforce that all numbers used in calculations are less
than $1.0$, then we know that the majority of the bits will produce
error that is also less than one, since 51 of the total 52 mantissa
bits will contribute an error less than $1.0$. We also see that some
of the exponent bits have the potential to contribute an error less than
$1.0$, which indicates if we can enforce or assume some properties of
the data used in the calculations e.g., data less than one, then we
can greatly increase the likelihood that a bit flip introduces an
error no greater than $1.0$. This phenomena is shown in
Figure~\ref{jje:fig:exp_bias_binary_perts}, where we show
empirically that numbers with exponent $2^0$ introduce a
small error  compared to the errors introduced with the exponent
$2^1$.

In summary, the exponent characterizes the error
introduced by the sign or mantissa should a bit flip occur.
As discussed in
\S~\ref{jje:sec:fault_model:characterization} and analytically
presented in Eq.~\eqref{jje:eq:magnitude_characterizes_mantissa}, we
are able to relate bit upsets to numerical error in terms of the
exponent of the original number.
Table~\ref{jje:table:bit_flips:scalar_error} summarizes these
analytical bounds for a bit upset in a scalar and
highlights the change in order of magnitude.
Now that we have characterized a fault in a scalar, we will present a
fault model centered around operations on scalars assuming one
will be perturbed.

\begin{table*}[htbp]\centering
\caption{Bit flip absolute error for a scalar  $\lambda$ represented using
IEEE-754 double precision, with
$\lambda~=~\lambda_\mathrm{exp}~\times~\lambda_\mathrm{frac}$. Where $\lambda_\mathrm{exp}$ is the
exponent $2^{x}$, and $\lambda_\mathrm{frac}$ is the fractional component.}
\label{jje:table:bit_flips:scalar_error}
\begin{tabular}{lcclc}
\toprule
Bit Location
& \multicolumn{2}{c}{Absolute Error}
					 & Bit Range
					 & $\Delta$ Order${}^\dagger$ \\
\cmidrule(r){2-3}
	& Scalar: $\Abs{\lambda - \tilde{\lambda}}$
	& Multiplication: $\Abs{\alpha\beta - \tilde{\alpha}\beta}$
	&
	& \\
\midrule
Mantissa
		& $\Abs{ \lambda_\mathrm{exp} (1 + 2^{j - 52} ) }$
		& $\Abs{ \alpha_\mathrm{exp} (1 + 2^{j - 52} )\beta }$
		& for $j = 0, \dotsc, 51$
		& $0$ \\
Exponent${}_{1 \to 0}$
		& $\Abs{ \lambda (1 - 2^{-2^j} )}$
		& $\Abs{ \alpha \beta (1 -2^{-2^j} ) }$
		& for $j = 0, \dotsc, 10$ and ${bit}_{j+52} = 1$
		& $-2^j$ \\
Exponent${}_{0 \to 1}$
		& $\Abs{ \lambda (1 - 2^{2^j} )}$
		& $\Abs{ \alpha \beta (1 - 2^{2^j} )}$
		& for $j = 0, \dotsc, 10$ and ${bit}_{j+52} = 0$
		& $+2^j$ \\
Sign
		& $\Abs{ 2\lambda  }$
		& $\Abs{ 2\alpha\beta  }$
		& 
		& $1$ \\ 
\bottomrule
\multicolumn{4}{l}{\footnotesize ${}^\dagger$ The change in order of magnitude.
$\phantom{\bigg(}$} 
\end{tabular}
\end{table*}

\subsection{Operation Centric Fault Model}
This work distinguishes itself from related work in the field of silent data
corruption by developing a fault model that is not based on perturbing arbitrary
memory locations.
We seek a fault model and experimental methodology that expresses all
possible errors, and not the expected error, which is what is obtained through
random sampling.
%

\subsubsection{Fault Model for Dot Products}
\label{jje:sec:fault_model:dotprod}
We now describe a realization of our fault model that describes the error that
could be injected if an operation in a dot product experiences a single bit
upset. We choose the dot product because it is a common operation, and because
we will use this model in \S~\ref{S:gmres} to model the
worst-case errors that could be injected into a phase of the GMRES algorithm.

Given two real-valued $n$-dimensional vectors $\Vect{a}, \Vect{b} \in \Rn$,
the dot product is defined as
\begin{equation}
c = \sum\limits_{i=1}^{n} c_i, ~~\text{where}~~c_i = a_i b_i.
\label{jje:eq:dot_prod}
\end{equation}
If we allow a single bit flip to impact the $i$-th element of the dot product,
then we have a perturbed solution $\tilde{c}$, which is the result of a
perturbation to either $a_i$, $b_i$, or $c_i$.
In the context of our fault model, this captures a bit upset impacting the
inputs to the multiplication operator, and it captures a bit upset in the
intermediate value, $c_i$, which is the input to the addition
operator.

Using Table~\ref{jje:table:bit_flips:scalar_error}, we have all of the
tools necessary to compose an absolute error model for a dot product,
i.e., addition is modeled by a fault in a scalar $\Abs{\alpha + \beta
- (\tilde{\alpha} + \beta)} = \Abs{\alpha - \tilde{\alpha}}$.
The potential change in order of magnitude is paramount. Consider an
exponent flip from $1 \to 0$.  These types of exponent bit flips
produce an error that is bounded above by the original magnitude of
the result, which can be viewed as ``zeroing out'' the term if a
perturbation occurs. Similar to a perturbed scalar, the mantissa can
contribute either no change in the order of magnitude, or in the worst
case a bit flip causes a carry, which will increment the order of
magnitude by one. The order of magnitude for a sign bit flip is
exactly the same as that of a perturbed scalar, which introduces
an error one order of magnitude larger than the result. These error
models can be thought of as the largest additive error that we can inject into a
dot product from a bit flip, e.g.,
\begin{equation}
\tilde{c} = \sum\limits_{i=1}^{n} a_i b_i + (\mathrm{error~term}).
\label{jje:eq:dot_prod_and_error}
\end{equation}

In summary, we have composed analytic models for the the absolute error that
could be introduced into a dot product. Our models are initially constructed
from the IEEE-754 Binary64 model, which we extended to express how a bit upset
impacts a singular double precision scalar. We then composed a model for the
multiplication operator, and analytically expressed the absolute error. Using
the absolute error, we have a model that explains \emph{how wrong} a dot product
can be, assuming a bit flip in one of the input vectors or in an intermediate
value.
Next, we refine these models to construct strict upper bounds on the error
introduced by a bit flip in a dot product.

\subsubsection{Error Bounds for a Bit Flip in a Dot Product}
\label{jje:sec:fault_model:dotprod_error}
The models presented in Table~\ref{jje:table:bit_flips:scalar_error}
make no assumptions about the bits present in the mantissa of the operands. This
is problematic if we want to consider \emph{all} possible errors that could be
introduced into a dot product. To account for the mantissa, and to
create strict upper bounds on the error, we will use the relation presented in
Eq.~\ref{jje:eq:magnitude_characterizes_mantissa}. From this relation, we know
that $\alpha \beta < 2^{\alpha_\mathrm{exponent} + 1} 
2^{\beta_\mathrm{exponent} + 1} $.  We can write this as
\begin{equation}
\alpha \beta  < 4 \alpha_\mathrm{exp} \beta_\mathrm{exp},
\label{jje:eq:fault_model:dotprod:mantissa_bound}
\end{equation}
where $\alpha_\mathrm{exp} = 2^{\alpha_\mathrm{exponent}}$. Using
Eq~\eqref{jje:eq:fault_model:dotprod:mantissa_bound}, we are able to account for the mantissa bits, but we can also show that a bit flip in
the sign is bounded by
Eq.~\eqref{jje:eq:fault_model:dotprod:mantissa_bound}. The sign bit introduces
an absolute error equivalent to incrementing the exponent of the
result
\begin{equation}
\alpha \beta  < 2 \alpha \beta < 4 \alpha_\mathrm{exp} \beta_\mathrm{exp},
\label{jje:eq:fault_model:dotprod:sign_bound}
\end{equation}
where $2 \alpha \beta$ is the potential error introduced should the sign bit be
perturbed, which must be smaller than the bound constructed for the mantissa.

By utilizing Eq.~\eqref{jje:eq:fault_model:dotprod:mantissa_bound}, we are able
to account for all possible mantissas and their potential faults, as well as a
perturbation to the sign bit. We will now discuss how to use this model to
understand the relationship between the data in an algorithm and the
distribution of potential errors that could occur should a bit flip in the
data.

\section{Fault Model Evaluation} \label{S:fault_model_eval}
In Section~\ref{S:fault_model} we proposed analytic models for errors
should a bit flip occur in IEEE-754 double precision data. We now
illustrate how data can impact the size of errors that a bit flip can
create. Consider the following sample vectors
\begin{alignat*}{2}
&
\Vect{u}_{\mathrm{small}} =
\left[\!\!\!\begin{array}{c}
0.5 \\
0.25
\end{array}\!\!\!
\right],~
&&
\Vect{v}_{\mathrm{small}} =
\left[\!\!\!\begin{array}{c}
0.25 \\
0.5
\end{array}\!\!\!
\right]
,~\text{and}~ \\
&
\Vect{u}_{\mathrm{large}} =
\left[\!\!\!\begin{array}{c}
2 \\
4
\end{array}\!\!\!
\right],~
&&
\Vect{v}_{\mathrm{large}} =
\left[\!\!\!\begin{array}{c}
4 \\
2
\end{array}\!\!\!
\right].
\end{alignat*}
If we compute the dot product $\lambda = \Vect{u}_{\mathrm{large}} \cdot
\Vect{v}_{\mathrm{large}}$, we have a finite number of potential errors should a
bit flip in the data of $\Vect{u}_{\mathrm{large}}, \Vect{v}_{\mathrm{large}}$,
or in an intermediate value in the summation. We can experience either
$\tilde{2} \times 4 + 4 \times 2$, $2 \times \tilde{4} + 4 \times 2$, or
$\tilde{8} + 8$. We have previously shown what $\tilde{2}$ can be (in
Figure~\ref{jje:fig:exp_bias_binary_perts}), but for clarity we will
state what the perturbed values could be (in
Figure~\ref{F:large_perturbed_values}).
\begin{figure}[htp]
\begin{equation*}
\tilde{2} = 
\left\{\!\!\!\begin{array}{c}
	2^{2}  \\ 
	2^{3}  \\
	2^{5}  \\
	2^{9}  \\
	2^{17} \\
	2^{33} \\
	2^{65} \\
	2^{129}\\
	2^{257}\\
	2^{513}\\
	\mathrm{Zero} \\
\end{array}\!\!\!\right\},\quad
\tilde{4} = 
\left\{\!\!\!\begin{array}{c}
	2^{1}  \\
	2^{4}  \\
	2^{6}  \\
	2^{10}  \\
	2^{18}  \\
	2^{34}  \\
	2^{66}  \\
	2^{130}  \\
	2^{258}  \\
	2^{514}  \\
	2^{-1020}  \\
\end{array}\!\!\!\right\},\quad
\tilde{8} = 
\left\{\!\!\!\begin{array}{c}
	2^{4}  \\
	2^{1}  \\
	2^{7}  \\
	2^{11}  \\
	2^{19}  \\
	2^{35}  \\
	2^{67}  \\
	2^{131}  \\
	2^{259}  \\
	2^{515}  \\
	2^{-1018}  \\
\end{array}\!\!\!\right\}
\end{equation*}
\caption{Example of perturbed values for large numbers.}
\label{F:large_perturbed_values}
\end{figure}
By inspection it is clear that substituting any of the above perturbed scalars
into the dot product will produce an absolute error greater than one in all cases,
and in the event one chooses to substitute the near zero perturbed
values, the absolute error of the dot product still has magnitude $8$, e.g., $\Abs{16 - (0 +
8)}$.

Alternatively, consider the vectors 
$\Vect{u}_{\mathrm{small}}$ and $\Vect{v}_{\mathrm{small}}$. If we compute
the dot product, $\lambda = \Vect{u}_{\mathrm{small}} \cdot
\Vect{v}_{\mathrm{small}} = 0.25$. Then we have possible values to
perturb: $\widetilde{0.5}$, $\widetilde{0.25}$, and
$\widetilde{0.125}$. We construct these from our model of a
perturbed scalar, and present the perturbed variants in
Figure~\ref{F:small_perturbed_values}.
\begin{figure}[htp]
\begin{equation*}
\widetilde{0.5} = 
\left\{\!\!\!\begin{array}{c}
	2^{0}  \\
	2^{-3}  \\
	2^{-5}  \\
	2^{-9}  \\
	2^{-17}  \\
	2^{-33}  \\
	2^{-65}  \\
	2^{-129}  \\
	2^{-257}  \\
	2^{-513}  \\
	2^{1022}  \\
\end{array}\!\!\!\right\},
\widetilde{0.25} = 
\left\{\!\!\!\begin{array}{c}
	2^{-3}  \\
	2^{0}  \\
	2^{-6}  \\
	2^{-10}  \\
	2^{-18}  \\
	2^{-34}  \\
	2^{-66}  \\
	2^{-130}  \\
	2^{-258}  \\
	2^{-514}  \\
	2^{1019}  \\
\end{array}\!\!\!\right\},
\widetilde{0.125} = 
\left\{\!\!\!\begin{array}{c}
	2^{-2}  \\
	2^{-1}  \\
	2^{-7}  \\
	2^{-11}  \\
	2^{-19}  \\
	2^{-35}  \\
	2^{-67}  \\
	2^{-131}  \\
	2^{-259}  \\
	2^{-515}  \\
	2^{1017}  \\
\end{array}\!\!\!\right\}
\end{equation*}
\caption{Example of perturbed values for small numbers.}
\label{F:small_perturbed_values}
\end{figure}
By inspection, $\widetilde{0.5}$ can contribute an absolute error to the dot
product larger than one only once, e.g., $\Abs{0.25 - (2^{1022}\times 0.25 +
0.125) }$. Likewise, $\widetilde{0.25}$ and $\widetilde{0.125}$ can perturb the
result of the dot product with error greater than one only once, and for all 3
cases the perturbation will change the result by hundreds of orders of
magnitude.

Returning to Figure~\ref{jje:fig:exp_bias_binary_perts} explains what causes bit
flips in the exponent to produce either a majority of large or small
errors. The binary pattern of the
stored biased exponent contains predominantly zeros for numbers greater than
one, and predominantly ones for numbers less than one.

One could also obtain primarily ones in the exponent as you approach
the extrema of the biased exponents, i.e., numbers larger than
$2^{512}$. In this case, the biased exponent does contain many ones, 
however, because the number is sufficiently large, i.e., the absolute
error will remain considerably large.  This is because if one ``zeroes out''
a perturbed element in the dot product, the error is proportional to
the magnitude of the result.



\subsection{Faults in the Mantissa or Sign}
The error generated by the mantissa or sign bits is relative to the
exponent of the number that the flip occurred in. If the exponent is
larger than one, then clearly the mantissa or sign bits can generate
an error larger than one.
Alternatively, if the values all are less than one, then mantissa
errors will produce errors less than one because $ 2^{-1} \times 1.x
\leq 2^0$. The errors from the sign bit cannot exceed $2$ since $ 2
\times 2^{-1} \times 1.x < 2^1$.

It is reasonable to consider that the mantissa generates a carry, as discussed
in \S~\ref{jje:sec:fault_model:dotprod_error}. To account for this we construct
a strict upper bound by incrementing the exponent of each element of the vectors
analyzed, similar to Eq.~\eqref{jje:eq:fault_model:dotprod:mantissa_bound}. For
example,
\begin{equation}
\Vect{u}_{\mathrm{original}} =
\left[\!\!\!\begin{array}{c}
2.12332 \\
1.24568
\end{array}\!\!\!
\right]~\Rightarrow~
\Vect{u}_{\mathrm{upper~bound}} =
\left[\!\!\!\begin{array}{c}
4 \\
2
\end{array}\!\!\!
\right].
\end{equation}
We then can evaluate our models on these vectors to determine a strict upper
bound on the errors we can experience in a dot product.

\subsection{Modeling Large Vectors}
\label{jje:sec:fault_model_eval:large_vectors}
We have shown how to exhaustively examine each element in a vector, and from
this analysis we can determine precisely which absolute errors we could
experience. Given large vectors, where the dimension $n$ may have millions or
billions of elements, exhaustively searching each element would be time
consuming, but it would also be a waste of time. As stated previously, there is a discrete
number of exponents supported by the IEEE-754 Binary64
specification. As we have previously shown, the exponent characterizes the
faults we can observe, so we only need to consider the 2046 possible biased
exponents and the special case of zero.
The perturbations that are possible can be determined independent of
concrete data values, e.g., we can precompute the perturbations and
absolute error because we know the relation stated in Eq.~\eqref{jje:eq:magnitude_characterizes_mantissa} and Eq.~\eqref{jje:eq:fault_model:dotprod:mantissa_bound}.

To analyze arbitrarily large vectors, we construct a lookup table for the 
absolute error in whatever operation we choose to model (we have chosen
products and addition). The table size is $2047 \times 2047$, and allows us to
consider the error introduced by performing an operation on two exponents,
which will map to a unique $ij$ location.

For example, consider the vectors
\begin{equation}
\Vect{u} =
\left[\!\!\!\begin{array}{c}
1.0 \\
1.2 \\
8.0 \\
0.125 \\
\end{array}\!\!\!
\right],~\text{and}~
\Vect{v} =
\left[\!\!\!\begin{array}{c}
0.125 \\
0.125001 \\
0.125002 \\
1.0 \\
\end{array}\!\!\!
\right].
\end{equation}
We first extract the biased exponents from the vectors
\begin{equation}
\Vect{u} \Rightarrow
\Vect{u}_\mathrm{exponent} =
\left[\!\!\!\begin{array}{c}
2^0 \times 1.0 \\
2^0 \times 1.x \\
2^3 \times 1.0 \\
2^{-3} \times 1.0 \\
\end{array}\!\!\!
\right] \Rightarrow
\Vect{u}_\mathrm{biased} =
\left[\!\!\!\begin{array}{c}
1023 \\
1023 \\
1026 \\
1020 \\
\end{array}\!\!\!
\right]
\end{equation}
Now, we determine an interval of possible values, and account for the
mantissa values that may have been truncated
\begin{equation}
u_i \in [1020, 1026] \subseteq [1020, 1027]~\text{for}~i=1,\dotsc,4.
\end{equation}

The range of biased exponents $[1020, 1027]$ will contain all possible values
that the original vector contained, and include one value that was larger than any in
the vector, the number corresponding to the biased exponent $1027$. Similarly,
we can compute the interval for $\Vect{v}$
\begin{equation}
\Vect{v} =
\left[\!\!\!\begin{array}{c}
0.125 \\
0.125001 \\
0.125002 \\
0.25 \\
\end{array}\!\!\!
\right] \!\!\Rightarrow\!\!
\left[\!\!\!\begin{array}{c}
2^{-3} \times 1.0 \\
2^{-3} \times 1.x \\
2^{-3} \times 1.x \\
2^{-2} \times 1.0 \\
\end{array}\!\!\!
\right] \!\!\Rightarrow\!\!
\left[\!\!\!\begin{array}{c}
1020 \\
1020 \\
1020 \\
1021 \\
\end{array}\!\!\!
\right],
\end{equation}
which leads to the interval we consider errors
\begin{equation}
v_i \in [1020, 1021] \subseteq [1020, 1022]~\text{for}~i=1,\dotsc,4.
\end{equation}
To allow us to analyze intervals efficiently, we create a
lookup table,
where each entry computes the relevant perturbations and absolute
errors for the operations being modeled. In the case of multiplication, the table has symmetry
because multiplication is commutative. In practice, computing the full table
($0,\dotsc,2046$) is simple and allows one to model errors for arbitrary
vectors.

A caveat of the above approach is that we must know the range of values that the
vector contains. This can be achieved by directly computing the min and max
values for each vector. Alternatively, an approximate range can be determined if
the ``length'' of the vector is known, e.g., the two-norm or \emph{if
we know that the data is normalized}, i.e., the two-norm is one.
One weakness to the proposed approach is that we do not consider a flip in the accumulating sum, which we
have left to future work. We also leave to future work analysis that
shows how many of these modeled errors lie within the rounding error
bound for pairwise sums.

\subsection{Summary}
We have shown that the range of values used in the dot product has a
direct impact on the size of the errors that can be observed. A general rule in
floating point algorithms has been to perform operations on numbers as close to
the same magnitude as possible, as doing so minimizes the loss of precision. We
have now shown that following this rule-of-thumb also gives the benefit of
making bit upsets generate a relatively small error when the numbers
are no larger than one.
 Next we present a motivating case study
that focuses exclusively on dot products, and then in
\S~\ref{S:gmres} we show how to exploiting data scaling in an
iterative solver.

\section{Case Study: Vector Dot Products} \label{S:case_study}
To begin our investigation, we assess the susceptibility of the
dot product of two $N$-dimensional vectors to a
silent bit flip in arithmetic. We make this choice since
many linear algebra operations can be decomposed into dot products,
for example, Gram-Schmidt orthogonalization or matrix-vector
multiplications.
\subsection{Computational Challenges}
Given a single double-precision number, there are 64 bits that may be
flipped. Extending this to an $N$-dimensional vector, we have $64N$
bits that are candidates for flipping. Accounting for different
numbers results in a very large sample space, and, therefore, we
utilized a hybrid CPU-GPU cluster and created a parallel code that
farms out specific Monte Carlo trials to various nodes. In this
context, a trial consists of creating two vectors, which is discussed
in the follow section, and then determining pass/fail statistics given
a bit flip on the input to the dot-product kernel. We utilized the
BLAS $ddot()$ routine, and aggregated the output for
post-processing in MATLAB. Ensuring a sampling error of less than
0.001, which is discussed in Section~\ref{section:dotprod:comparison},
required nearly 400,000 CPU hours parallelized over the
processors of a 1700 core cluster of AMD 6128 Opterons.
The large
search space coupled with ensuring statistically significant results
highlights why an analytic approach is not only more efficient, but
may be required for more advanced methods and data structures, e.g.,
matrices and linear solvers.
\subsection{Monte Carlo Sampling}
We next develop a better understanding of how vector magnitudes impact
the expected absolute error should a bit perturb a dot product. To
conduct Monte Carlo sampling, we must first determine a mechanism for
tallying success, and we must define \textit{success and failure}.
\begin{itemize}
\item \textbf{Vector Creation} \\
  \hspace*{-1em}1)~Mantissa generated randomly using C stdlib rand(). \\
  \hspace*{-1em}2)~For each vector, we fix each element's magnitude to
  the bit pattern $2^{-50}$ to $2^{50}$ (101 bit patterns). This
  corresponds to the base ten numbers in the range $8.8\times
  10^{-16}$ to $1.1\times 10^{+15}$. This range was chosen because $2^{-50}$
  roughly is the machine precision. The numbers
  in this range are utilizing the highest \emph{precision} that
  Binary64 offers.
\item \textbf{Sample definition and Error Calculation} \\
  \hspace*{-1em}1)~A random sample is defined by generating two random
  $N$ length vectors and computing the absolute error considering all
  possible $2\times 64 \times N$ bit flips. \\
  \hspace*{-1em}2)~A tally is defined by failure, which we define to
  be any absolute error that is greater than $1$. \\
  \hspace*{-1em}3)~An empirical estimate of the expected absolute
  error is computed by dividing the number of failures by the number
  of bits considered times the vector length times 2 times the number
  of random samples ($M$) taken for a given magnitude combination, i.e.,
  $failures / (2 \times 64 \times N \times M)$.
\item \textbf{Visualization} \\
  \hspace*{-1em}1)~To visualize the expected absolute error, we
  construct tallies for each magnitude combination, i.e., $101 \times 101$
  unique combinations, and each combination is sampled $M$ times. \\
  \hspace*{-1em}2)~We summarize this information in a surface plot,
  where the x- and y-axes denote the $\log_2$ of the relative
  magnitude of the vector \Vect{u} and \Vect{v}, respectively. The
  height of the surface plot indicates the probability of seeing an
  absolute error larger than $1$.
\end{itemize}
\begin{figure*}[tpb]
\centering
\begin{subfigure}[b]{0.34\textwidth}
    \centering
    \includegraphics[width=\textwidth]{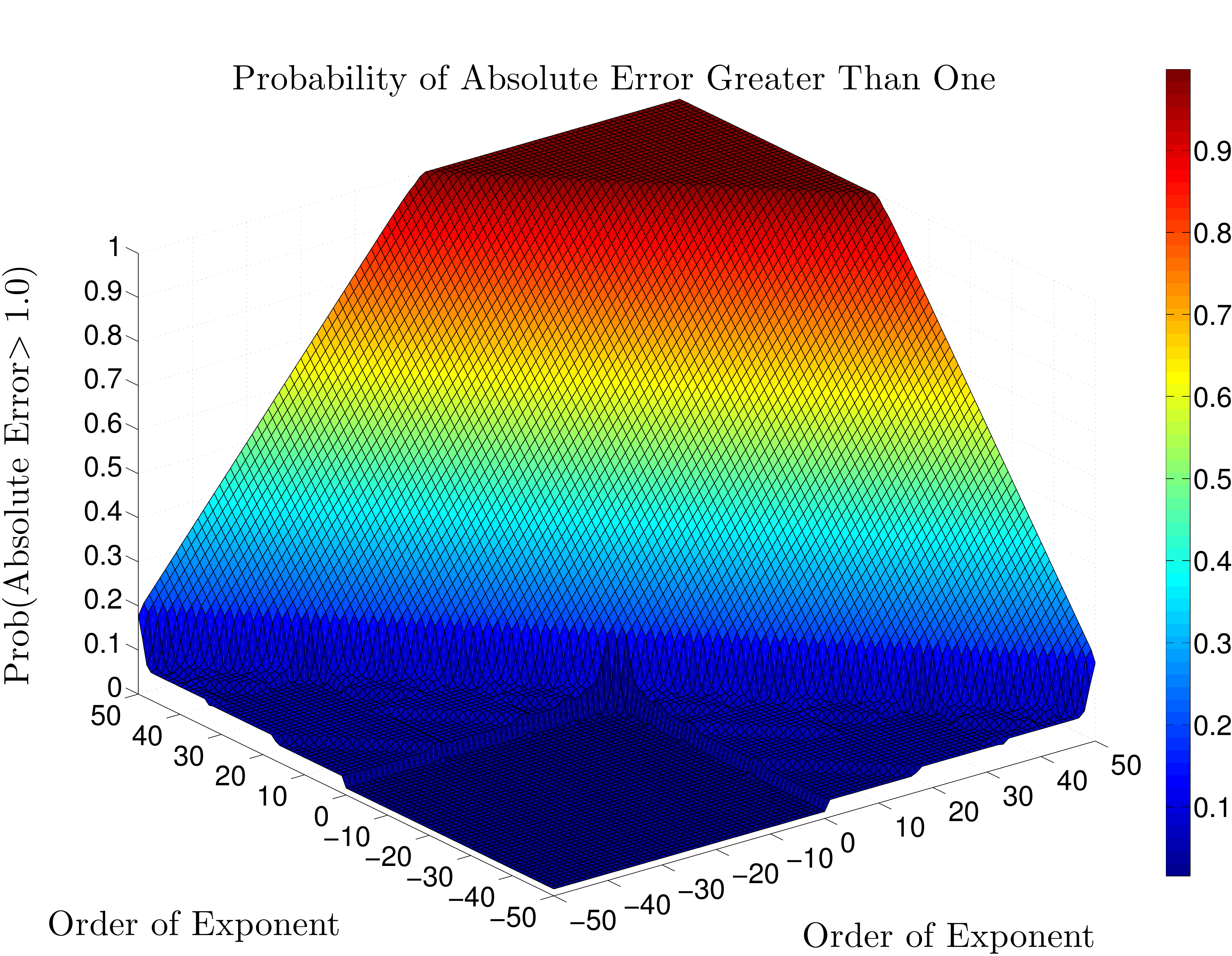}
    \caption{Monte Carlo experiment computing dot products with
    vectors of various magnitudes. Failure is defined to be an
    absolute error larger than 1.}
    \label{F:surfs:10000}
\end{subfigure}
\hfil
\begin{subfigure}[b]{0.28\textwidth}
    \centering
    \includegraphics[width=\textwidth]{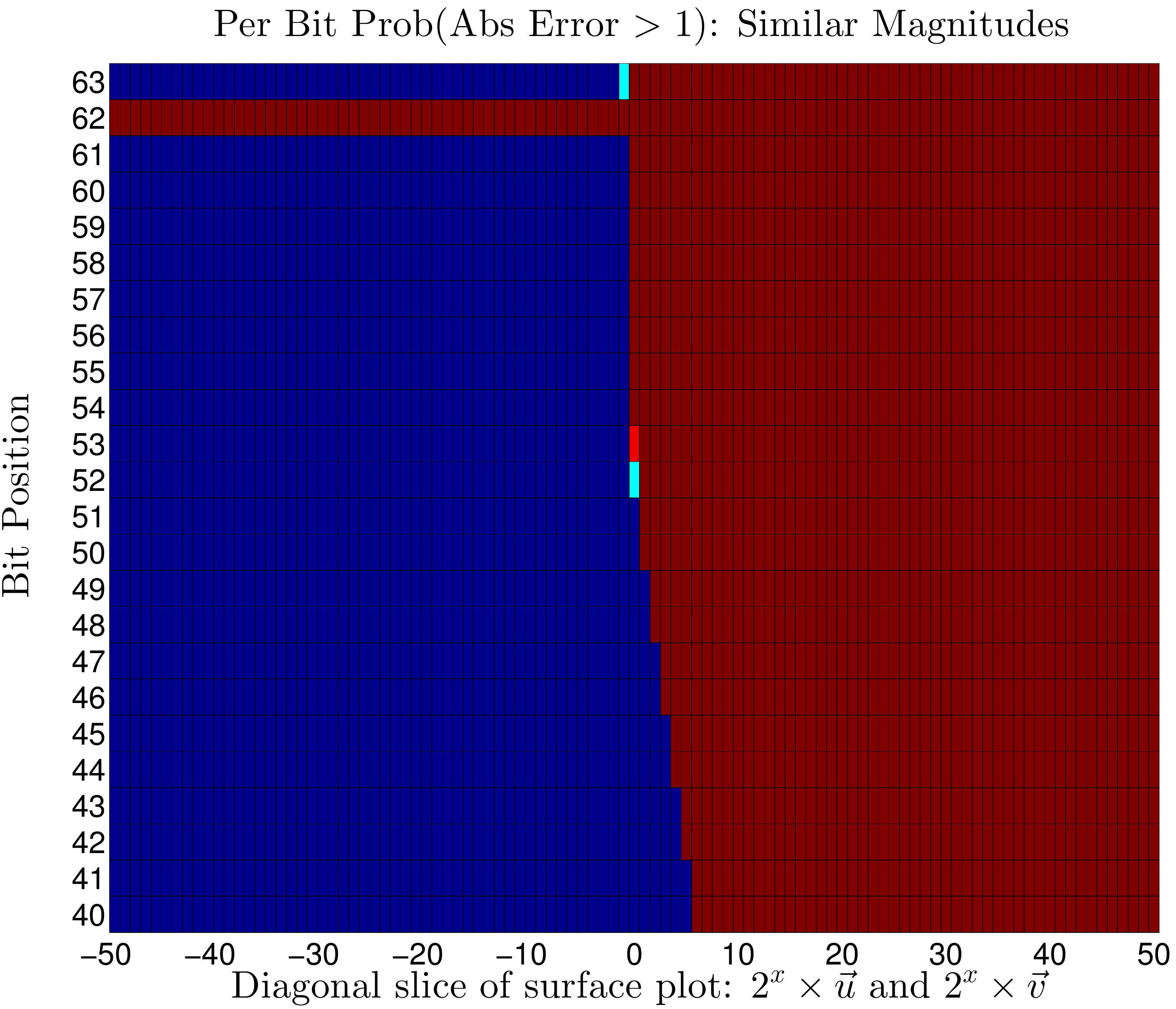}
    \caption{Dot product with vectors containing similar
    magnitudes.\\ \\}
    \label{F:slice:same:10000}
\end{subfigure}
\hfil
\begin{subfigure}[b]{0.28\textwidth}
    \centering
    \includegraphics[width=\textwidth]{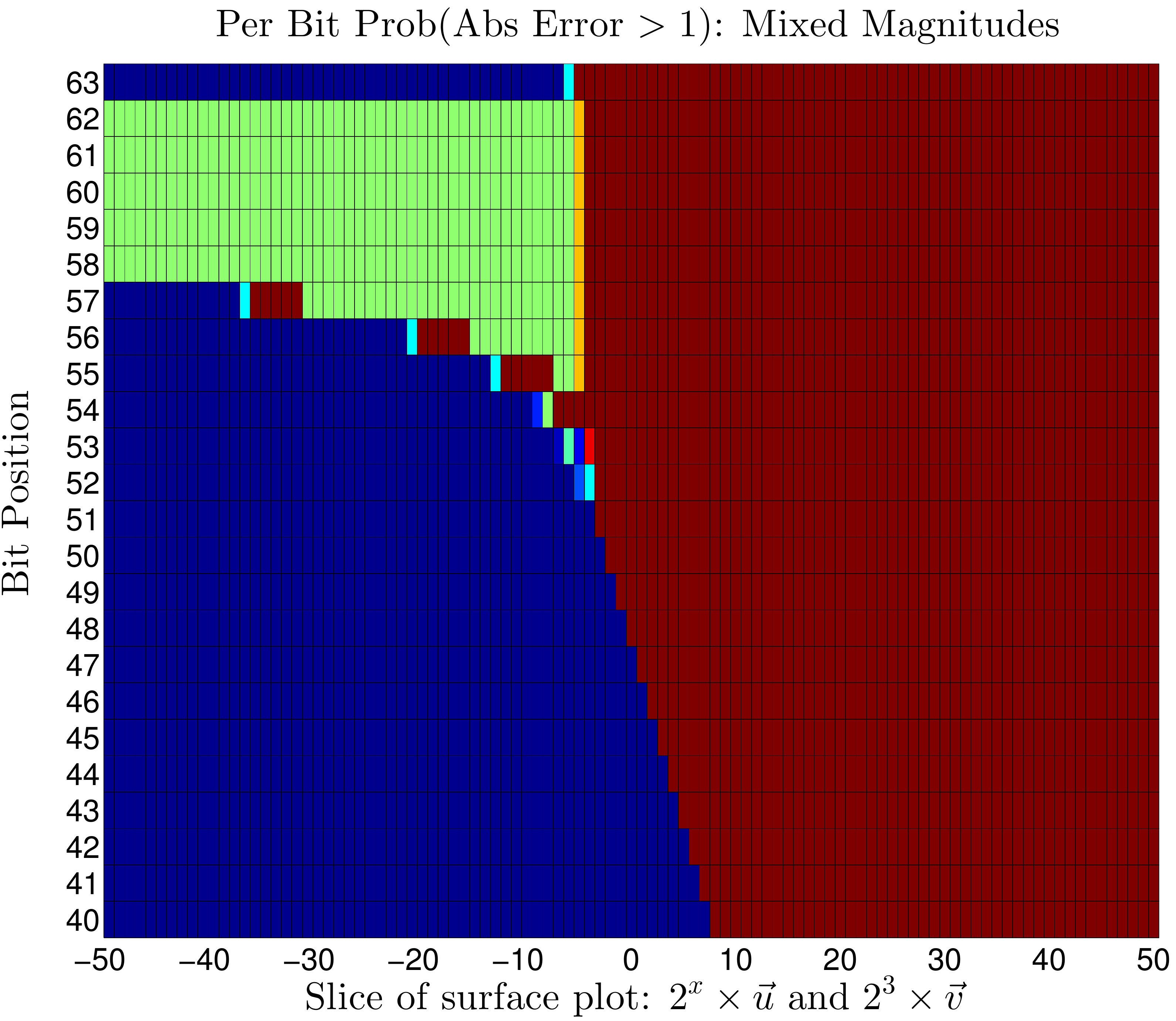}
    \caption{Dot product containing mixed magnitudes. \\ \\}
    \label{F:slice:diff:10000}
\end{subfigure}
\caption{Probability of observing an absolute error larger than 1.}
\label{F:surfs}
\end{figure*}
Figure~\ref{F:surfs:10000} presents a surface plot as described in the
Visualization bullet. To interpret this graph, the x-axis indicates
the magnitude that all elements of the vector $\Vect{u}$ were forced to
have while the mantissa was randomly generated. Likewise, the y-axis
indicates the magnitude that all elements of the vector $\Vect{v}$
were forced to have. Each x-y intersection represents 1,000,000 random
vector samples, where the dot product was computed and failures
tallied. The height of the surface at an $(x,y)$ location indicates
the probability of observing an absolute error larger than
$1$ given a single bit flip.  From this
surface, one may immediately recognize the unusual structure of
these graphs: When both vectors have magnitudes larger than $2^0$, the
probability of failure is noticeably higher; yet, when both vectors
have magnitudes less than or euqal to $2^0$, the probability of
failure is approaching zero.

The key finding presented in Figure~\ref{F:surfs:10000}, is that when
we operate on vectors that are normalized, e.g., values in the range
$[0,1]$, we have a very low probability of seeing a large error
should a bit flip occur. The lowest probability, i.e., the
flat region in the quadrant $[0,-50] \times [0,-50]$, is precisely
$\mathrm{Prob}(\mathrm{Abs~Error} > 1) = 0.015625$, which is $1/64$.
The single bit that can introduce absolute error larger than one is
the most significant exponent bit. Also, should the most significant
exponent bit flip, the error is quite large and can be detected
\cite{elliott13gmres}.

\subsection{Per Bit Analysis of Surface Plot}
\label{section:dotprod:per_bit_analysis}

To better understand the structure of the surface plot,  we take two
slices of the surface and look at the per-bit probability of a failure
(Figures~\ref{F:slice:same:10000}~and~\ref{F:slice:diff:10000}).
The slices chosen feature dot products of vectors with
\textit{similar} relative magnitudes and dot products
of vectors of many magnitudes (the x-axis) with a vector that contains
magnitudes up to $2^3$.
Intuitively, these figures slice from the back-most corner of
Figure~\ref{F:surfs:10000}  to the front for similar magnitudes
(Figure~\ref{F:slice:same:10000}), and they slice from the left to
right for Figure~\ref{F:slice:diff:10000}. 

We have shown why this shape should be expected in
Figure~\ref{jje:fig:exp_bias_binary_perts}, and in the example
presented in \S~\ref{S:fault_model_eval}.
This feature is an artifact of how the exponent is implemented in the
IEEE-754 specification, i.e., a biased exponent. The lowest
probability presented in the surface plot is $0.015625 = 1/64$, we
can graphically show this in Figure~\ref{F:slice:same:10000}, where one
can see that bit  \#62 (2\textsuperscript{nd} from the top), is the
only bit that can contribute large error. We also show that the sign
and mantissa bits can not introduce large error when values are in the
range $[0,1]$.

Conversely, Figure~\ref{F:slice:diff:10000} shows that
when mixing large and small values, we expect to see large errors for
faults. The green shading in
Figure~\ref{F:slice:diff:10000} (upper left quadrant) indicates a
roughly 50\% chance that we see an absolute error larger than one. The
reason for this is that values larger than 2 have a binary pattern
that introduces large error most of the time (recall
Figure~\ref{F:large_perturbed_values}).The increased likelihood of
large error from the large numbers, coupled with the low chance from small
numbers, creates a scenario where it is equally likely to see both
large and small absolute errors. The more we deviate from operating on
numbers in the range $[0,1]$, the closer we get to having a 50/50
chance of seeing a large error (see the mantissa bits slowly becoming
green as well).

\subsection{Comparison of the Analytic Model and Monte Carlo Sampling}
\label{section:dotprod:comparison}
\begin{figure}[!tb]%
\centering
\includegraphics[width=\columnwidth]{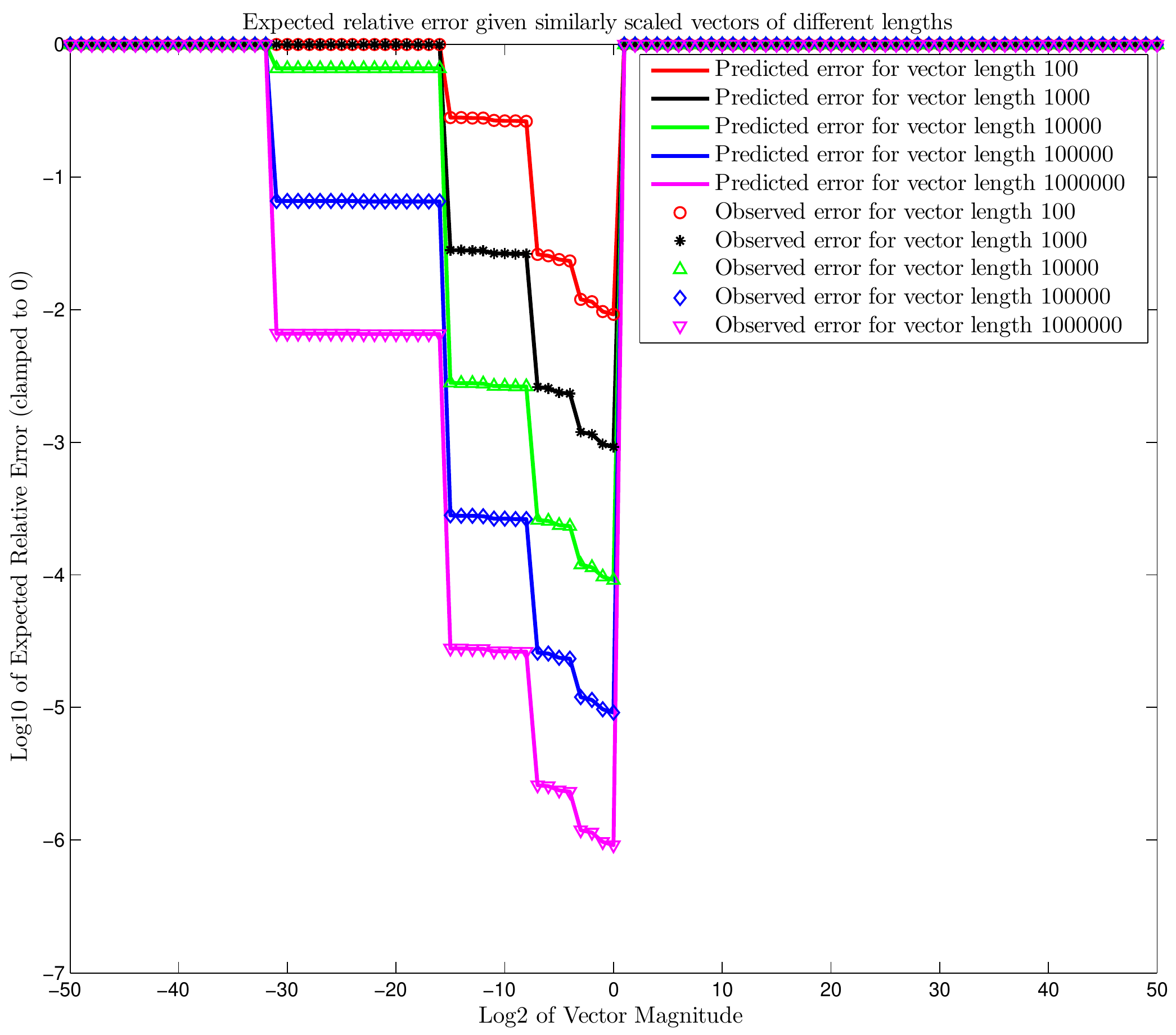}%
\caption{Comparison of observed error caused by a flip in the
  exponent, excluding the most significant bit, for sampled vector
  sizes having similar relative magnitudes.}%
\label{figure:analytical:comp}%
\end{figure}

In Figure~\ref{figure:analytical:comp}, we compare the error observed
while performing Monte Carlo sampling with the expected error
computed from our model. 
We sampled up to $M=$ 1 million random vectors per data point, which
implies a Monte Carlo error of $error_{MC} = 1/\sqrt{M} \approx
0.001$. We observe a perfect fit, which is to be expected because we
have analytically shown that the exponent bits dictate the size of the
absolute error we will observe. Even with random sign and mantissa
bits evaluated, we see that the likelihood of experiencing a large
error is entirely determined by the exponent bits.

\section{Extension to Matrices and Iterative Solvers} \label{S:gmres}
Having recognized that dot products on numbers less than one can
produce errors less than one, we will relate this idea to matrix
equilibration. We then provide an example of how to use this concept
in an sparse iterative solver (GMRES), while exhaustively counting the
possible errors that can be introduced.

\subsection{Matrix Equilibration}
The idea of scaled vectors is analogous to vector normalization, i.e.,
$\TwoNorm{\Vect{u}} = 1$. Applied to matrices in the context of
solving linear systems, scaling takes the form of \emph{matrix
  equilibration}: for a matrix $\Mat{A}$, scale the rows and columns
such that $\InfNorm{\Mat{A}} = 1$.  Scaling can also be performed
before a matrix is created, for example the equations leading to the
matrix can be scaled prior to assembling a matrix.  To scale a sparse
matrix after its creation, we use a sparse matrix implementation of
LAPACK's equilibration routine DGEEQU \cite{jje:lapackUsersGuide}.
Equilibration does not cause fill, i.e., it will not increase the
number of non-zeros. In general, equilibrating a matrix is only
beneficial, but equilibration may not be practical in all cases.

\subsection{GMRES}

The Generalized Minimum Residual method (GMRES) of Saad and Schultz
\cite{jje:saad1986gmres} is a Krylov subspace method for solving
large, sparse, possibly non-symmetric linear systems
$\Mat{A}\Vect{x}=\Vect{b}$.  GMRES is based on the Arnoldi process
\cite{mfh:arnoldi1951principle}, which uses orthogonal projections and
basis vectors normalized to length one. Arnoldi and GMRES relate to
this work because the orthogonalization phase of Arnoldi is often
Modified Gram-Schmidt or Classical Gram-Schmidt, which are dot product
heavy kernels.


We present the GMRES algorithm in Algorithm~\ref{jje:alg:gmres}.
The Arnoldi process is expressed on
Lines~\ref{jje:alg:gmres:arnoldi_start}--\ref{jje:alg:gmres:arnoldi_end}
in Algorithm~\ref{jje:alg:gmres}.  At its core is the Modified
Gram-Schmidt (MGS) process, which constructs a vector orthogonal to
all previous basis vectors $\Vect{q}_i$. The MGS process begins on
Line~\ref{jje:alg:gmres:mgs_start} and completes on
Line~\ref{jje:alg:gmres:mgs_end}. We now describe how we instrument
the orthogonalization phase and count the absolute errors that
\emph{could} be injected.

\begin{algorithm}
\caption{GMRES}
\label{jje:alg:gmres}
\begin{algorithmic}[1]
\Input{Linear system $\Mat{A}\Vect{x}=\Vect{b}$ and initial guess $\Vect{x}_0$}
\Output{Approximate solution $\Vect{x}_m$ for some $m \geq 0$}
\State{$\Vect{r}_0 := \Vect{b} - \Mat{A} \Vect{x}_0$}\Comment{Initial residual vector}
\State{$\beta := \TwoNorm{\Vect{r}_0}$, $\Vect{q}_1 := \Vect{r}_0 / \beta$}
\For{$j = 1, 2, \dots$ until convergence}
\label{jje:alg:gmres:arnoldi_start}
  \State{$\Vect{v}_{j+1} := \Mat{A} \Vect{q}_j$}\Comment{Apply the matrix $A$}
  \label{jje:alg:gmres:orthog_start_vec}
  \For{$i = 1, 2, \dots, j$}\Comment{Orthogonalize}
  \label{jje:alg:gmres:mgs_start}
    \State{$h_{i,j} := \Vect{q}_i \cdot \Vect{v}_{j+1}$}
    \label{jje:alg:gmres:hij}
    \State{$\Vect{v}_{j+1} := \Vect{v}_{j+1} - h_{i,j} \Vect{q}_i $}
    \label{jje:alg:gmres:wi}
  \EndFor
  \label{jje:alg:gmres:mgs_end}
  \State{$h_{j+1,j} := \TwoNorm{\Vect{v}_{j+1}}$}
  \label{jje:alg:gmres:hjp1}
  \If{$h_{j+1,j} \approx 0$} \label{jje:alg:gmres:happyCheck}
       \State{Solution is
       $\Vect{x}_{j-1}$}\label{jje:alg:GMRES:happy_breakdown}\Comment{Happy
       breakdown}
       \State{\textbf{return}}
  \EndIf
  \State{$\Vect{q}_{j+1} := \Vect{v}_{j+1} / h_{j+1,j}$}\Comment{New basis
  vector}
  \label{jje:alg:gmres:arnoldi_end}
  \State{$\Vect{y}_j := \argminy \TwoNorm{ \Mat{H}(\IndexRange{1}{j+1},
  \IndexRange{1}{j}) \Vect{y} - \beta
  \Vect{e}_1}$}
  \State{$\Vect{x}_j := \Vect{x}_0 + [\Vect{q}_1, \Vect{q}_2, \dots, \Vect{q}_j]
  \Vect{y}_j$}\Comment{Compute solution update}
\EndFor
\end{algorithmic}
\end{algorithm}

\subsection{Instrumentation and Evaluation}\label{SS:Classes}
To demonstrate the benefit of data scaling we have chosen 3 test
matrices.
We instrument the code and for each dot product in the
orthogonalization phase we determine an interval that describes the
range of values possible in the vectors. Then using our fault
model, we compute the absolute errors that are possible.
 Since we know the basis vectors
($\Vect{q}_i$) are normal, the intervals for the values in the vectors
are [0,1]. We compute the min and max for the unknown vector
$\Vect{v}$, and this determines the interval for the values in
$\Vect{v}$.
We use the intervals and our fault model to evaluate all absolute
errors that can be introduced from a single bit flip in the input
vectors. We classify the absolute error into four classes:
\begin{enumerate}
  \item Absolute error less than 1.0,
  \label{jje:sec:gmres:eval:err:less1}
  \item Absolute error greater than or equal to 1.0, but less than or equal to
  $\TwoNorm{\Mat{A}}$,
  \label{jje:sec:gmres:eval:err:grey}
  \item Absolute error greater $\TwoNorm{\Mat{A}}$.
  \label{jje:sec:gmres:eval:err:detectable}
  \item Error that is non-numeric, e.g., Inf or NaN.
  \label{jje:sec:gmres:eval:err:nonnumeric}
\end{enumerate}
We choose to include the 2\textsuperscript{nd} class of errors due to
recent work by by Elliott et al. \cite{elliott13gmres} that
demonstrates how to use a norm bound on the Arnoldi process to filter
out large errors in orthogonalization.

Classes~\ref{jje:sec:gmres:eval:err:less1}~and~\ref{jje:sec:gmres:eval:err:grey}
are \emph{undetectable}, while
Classes~\ref{jje:sec:gmres:eval:err:detectable}~and~\ref{jje:sec:gmres:eval:err:nonnumeric}
are detectable. Our goal is to ensure that should a bit flip, the error falls
into Classes~\ref{jje:sec:gmres:eval:err:less1},~\ref{jje:sec:gmres:eval:err:detectable},~and~\ref{jje:sec:gmres:eval:err:nonnumeric}
while minimizing or eliminating the occurrence of
Class~\ref{jje:sec:gmres:eval:err:grey} errors. We refer to
Class~\ref{jje:sec:gmres:eval:err:grey} errors as the \emph{grey area}, as they
are undetectable errors that we consider to be large.

\subsubsection{Sample Problems}

We have chosen three sample matrices to demonstrate our technique.  To
ensure reproducibility, we did not create any of these matrices from
scratch, rather we used readily available matrices. The first matrix
arises from a second-order centered finite difference discretization
of the Poisson equation. We generated this matrix using MATLAB's
built-in Gallery functionality.
The second matrix, CoupCons3D, presents a more realistic linear
system.  It comes from the University of Florida Sparse Matrix
Collection \cite{jje:davis2011university} and arises from a fully
coupled poroelastic problem.  The matrix is symmetric in pattern, but
not symmetric in values.  It is also fairly large, and has explicitly
stored zero values. The matrix is poorly scaled, with a mixture of
large and small values.
The final matrix, mult\_dcop\_03, is also from the Florida Sparse
Matrix Collection.  It arises from a circuit simulation problem, and
has good scaling inherently.
We have
summarized the characteristics of each matrix in Table~\ref{jje:table:gmres:sample_matrices}.
\begin{table}[htp]\centering
\caption[labelInTOC]{Sample Matrices}
\label{jje:table:gmres:sample_matrices}
\begin{tabular}{lrrr}
\toprule
Properties & Poisson100 & CoupCons3D & mult\_dcop\_03 \\
\midrule
number of rows  &
10,000 &
416,800 &
25,187 \\
number of columns  &
10,000 &
416,800 &
25,187 \\
nonzeros  &
49,600 &
17,277,420 &
193,216 \\
structural full rank  & 
yes &
yes & 
yes \\
explicit zero entries  &
0 &
5,044,916 & 
0 \\
type  &
real &
real &
real \\
structure  &
symmetric &
nonsymmetric &
nonsymmetric\\
positive definite  & 
yes &
no &
no \\
\bottomrule
\end{tabular}
\end{table}

We now scale the Poisson and CoupCons3D matrices and right-hand side
vectors such that they are equilibrated.
Table~\ref{jje:table:gmres:sample_matrices:norms} summarizes the norms
for each of our test matrices. We use the infinity norm ($\InfNorm{A}
\approx 1$) to measure whether a matrix is well scaled.  One can see
that both the Poisson and mult\_dcop\_03 matrices have infinity norms
not too much larger than one, while the CoupCons3D matrix is
inherently poorly scaled.  Our equilibration code ran out of memory
when attempting to equilibrate mult\_dcop\_03, but it is already well scaled.
\begin{table}[htp]\centering
\caption{Norms of Sample Matrices ${}^\dagger$}
\label{jje:table:gmres:sample_matrices:norms}
\begin{tabular}{lrrrr}
\toprule
Norm & \multicolumn{2}{c}{Poisson Equation} & \multicolumn{2}{c}{CoupCons3D} \\
\midrule
& No Scaling & Scaling & No Scaling & Scaling \\  
$\InfNorm{\Mat{A}}$
& $8.0$ & $2.0$ & $1.30 \times 10^6$ & $1.0$ \\
$\TwoNorm{\Mat{A}}$
& $7.999$ & $1.999$ & $1.20 \times 10^6$ & $1.0$ \\
$\FNorm{\Mat{A}}$
& $4.46 \times 10^2$ & $1.12 \times 10^2$ & $2.75 \times 10^6$ & $2.91 \times
10^2$ \\
\bottomrule
\multicolumn{5}{l}{\footnotesize ${}^\dagger$ mult\_dcop\_03 has
$\InfNorm{\Mat{A}} = 35.5$.$\phantom{\bigg(}$}
\end{tabular}
\end{table}

\subsection{Results}

\begin{figure*}[htb]
\centering
\begin{subfigure}[b]{0.3\textwidth}
    \centering
    \includegraphics[width=.75\textwidth]{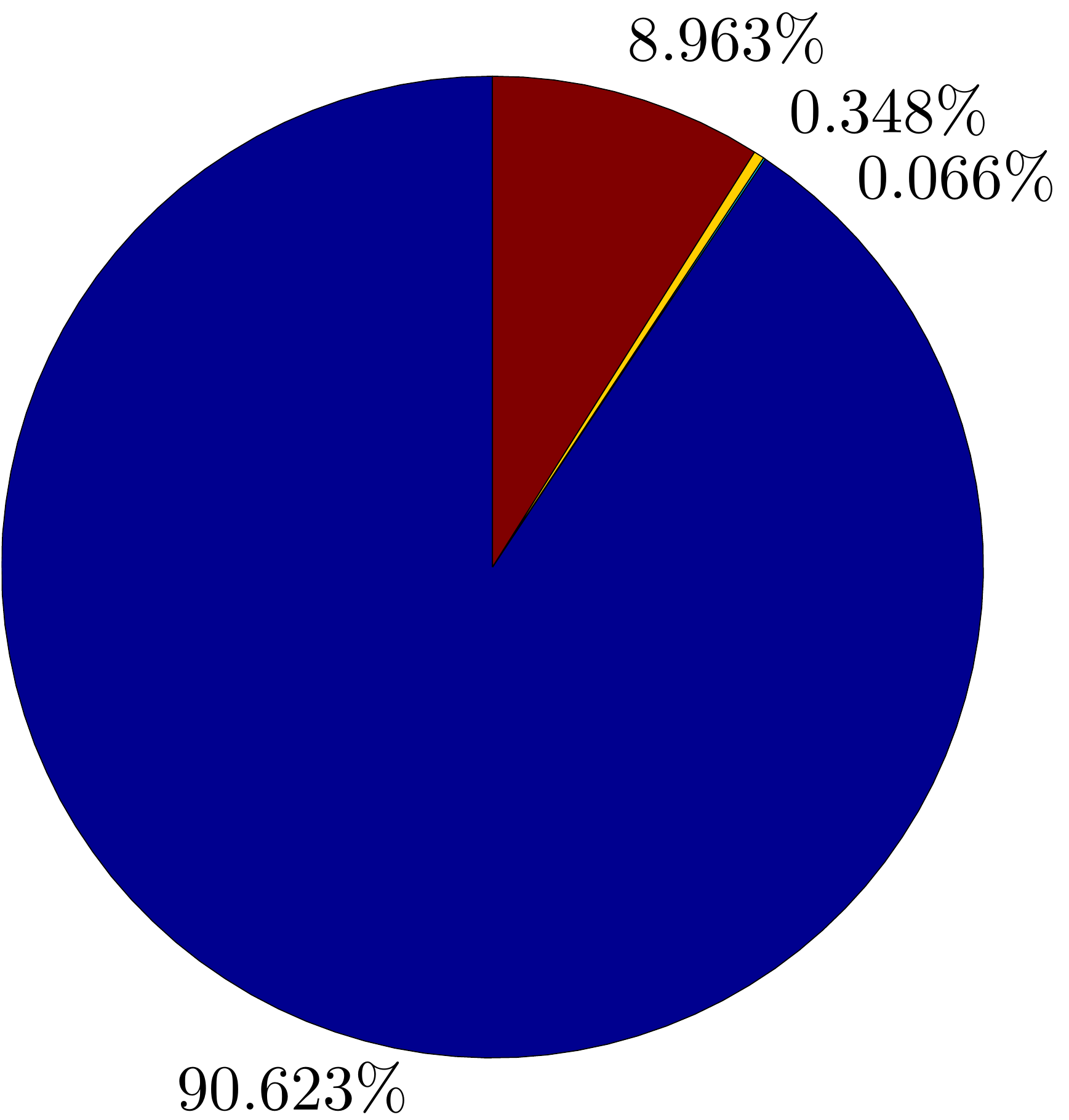}
    \caption{Poisson, no scaling.}
    \label{F:pie:noscaling:poisson}
\end{subfigure}
\hfil
\begin{subfigure}[b]{0.3\textwidth}
    \centering
    \includegraphics[width=.75\textwidth]{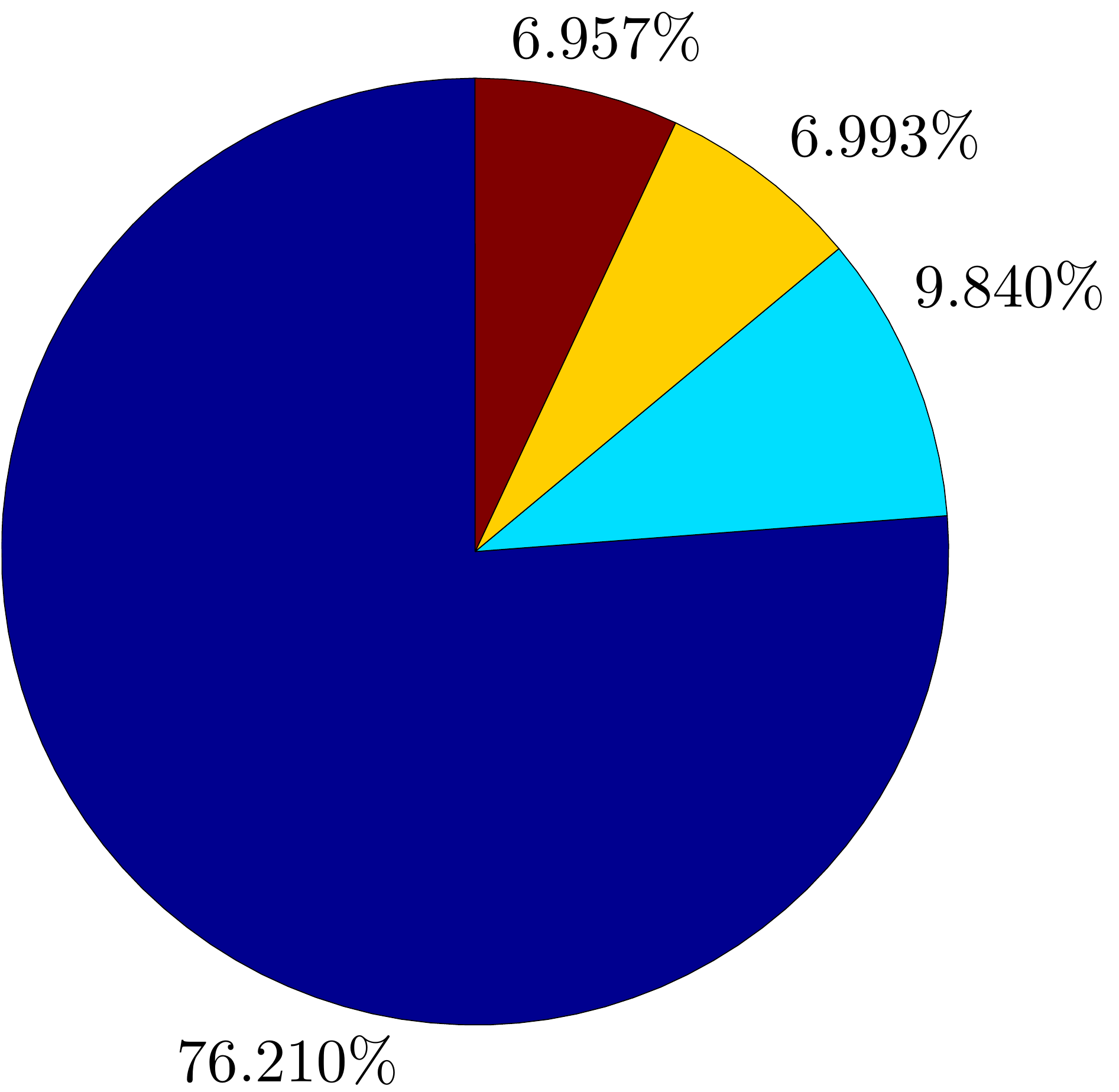}
    \caption{CoupCons3D, no scaling.}
    \label{F:pie:noscaling:coup}
\end{subfigure}
\hfil
\begin{subfigure}[b]{0.3\textwidth}
    \centering
    \includegraphics[width=\textwidth]{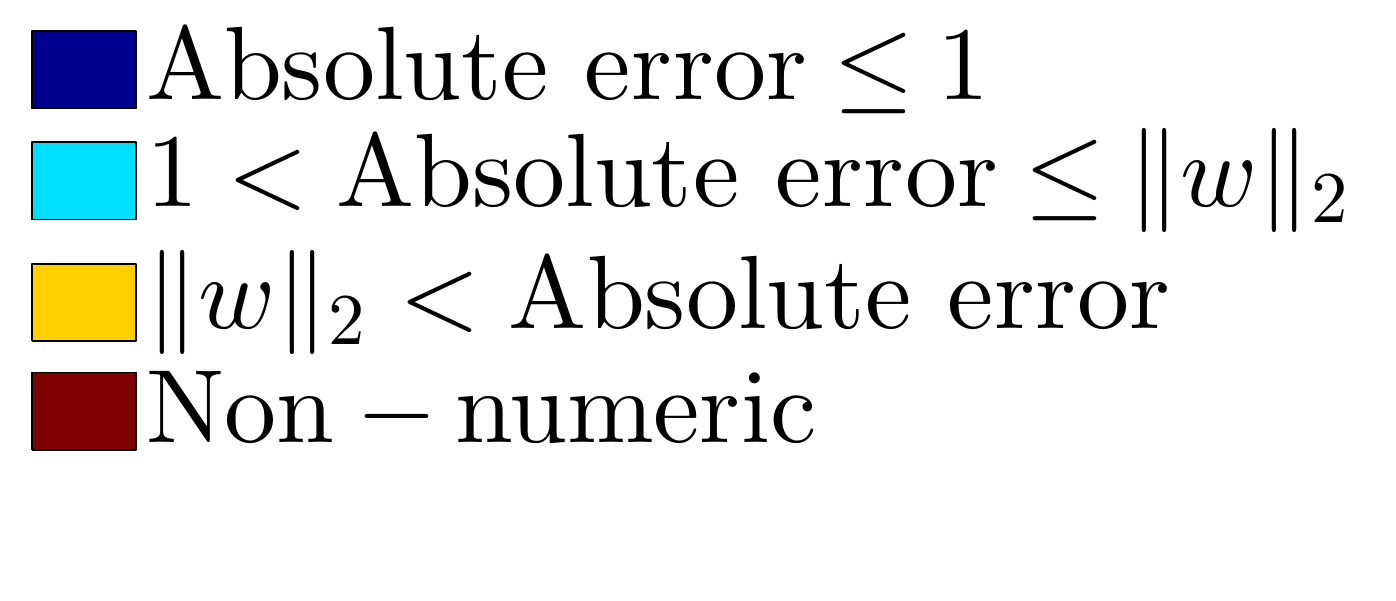}
    \vspace*{.25in}
\end{subfigure}
\\
\begin{subfigure}[b]{0.3\textwidth}
    \centering
    \includegraphics[width=.75\textwidth]{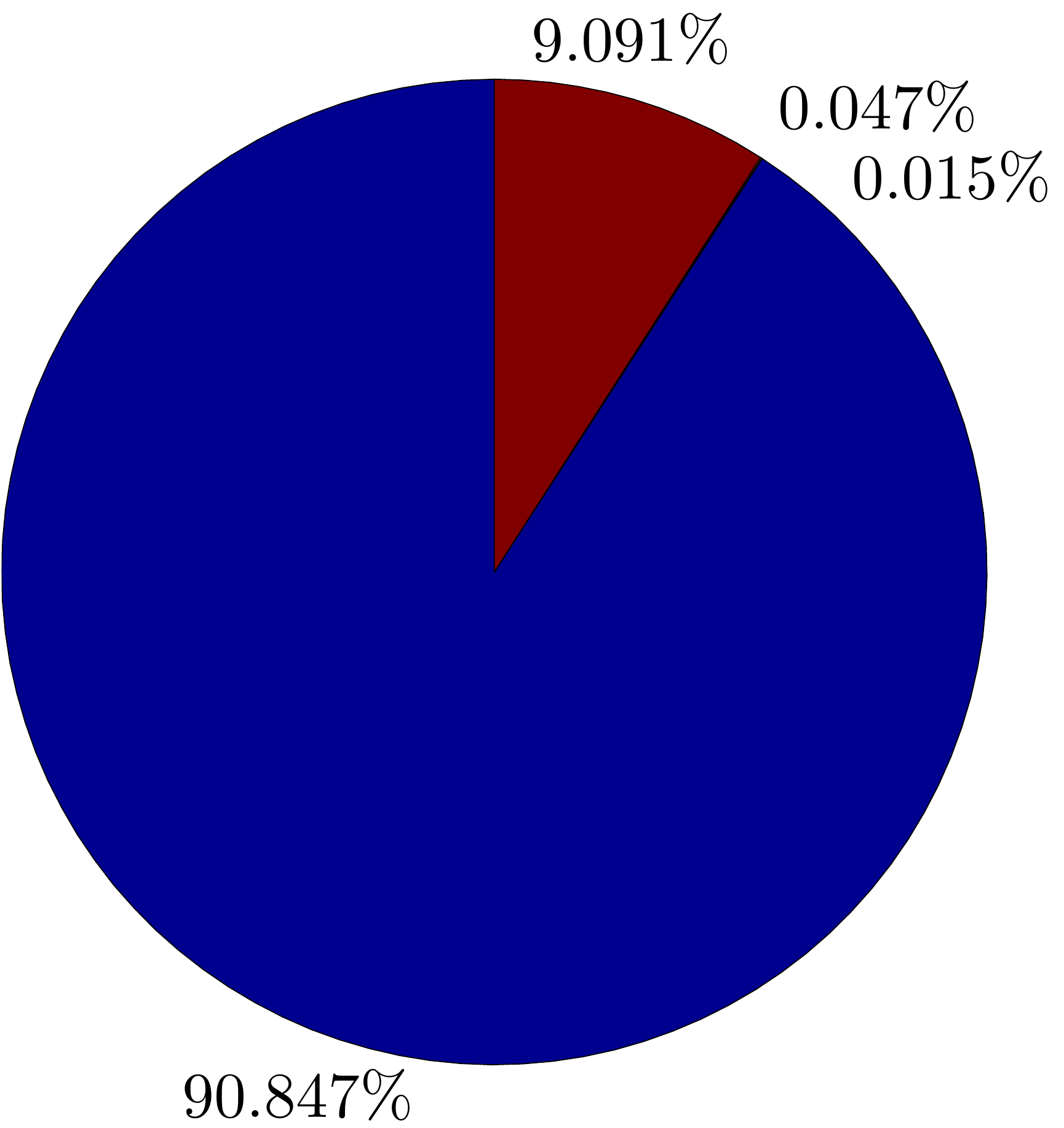}
    \caption{Poisson, equilibrated}
    \label{F:pie:noscaling:poisson:equib}
\end{subfigure}
\hfil
\begin{subfigure}[b]{0.3\textwidth}
    \centering
    \includegraphics[width=.75\textwidth]{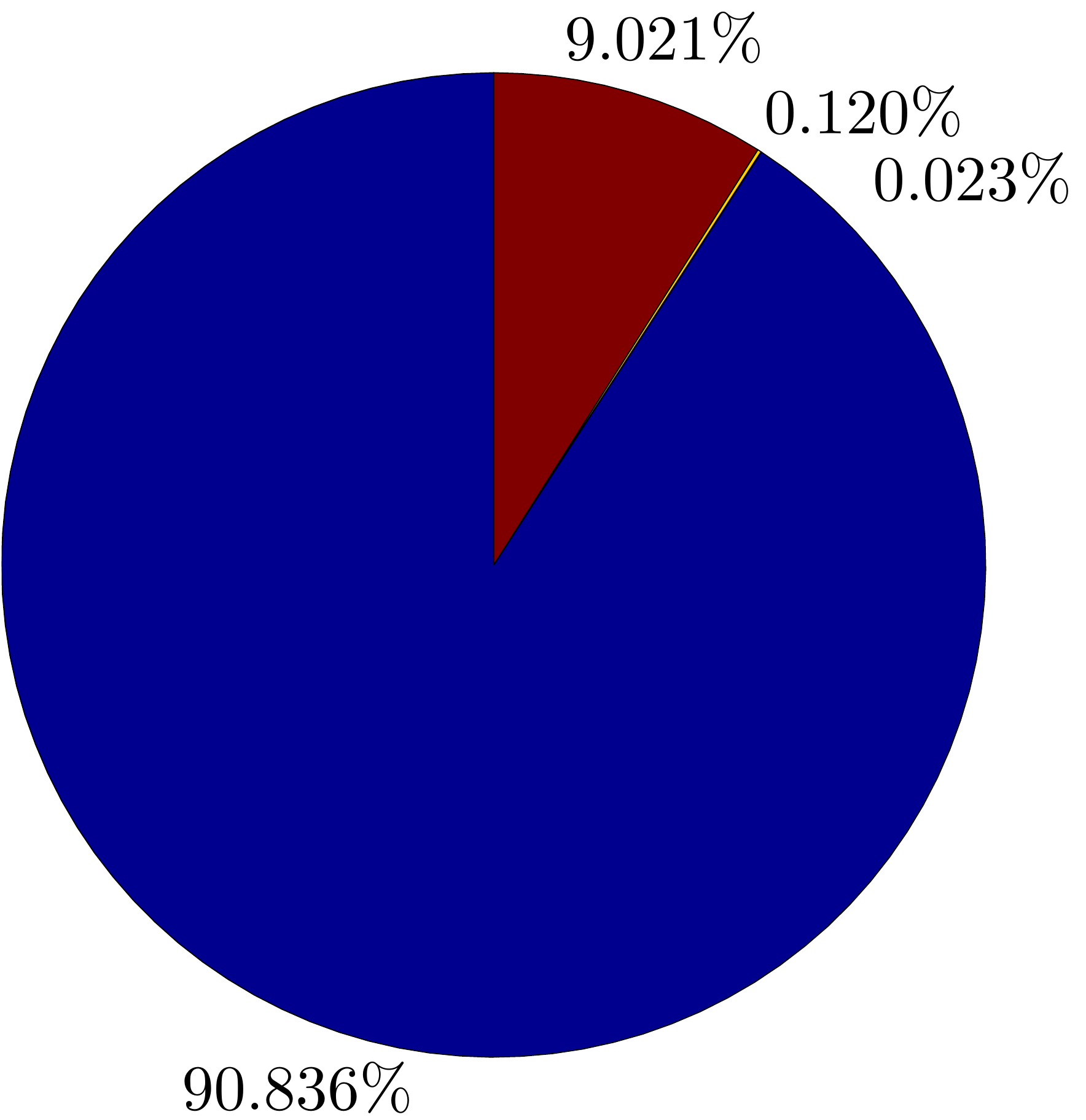}
    \caption{CoupCons3D, equilibrated.}
    \label{F:pie:noscaling:coup:equib}
\end{subfigure}
\hfil
\begin{subfigure}[b]{0.3\textwidth}
    \centering
    \includegraphics[width=.75\textwidth]{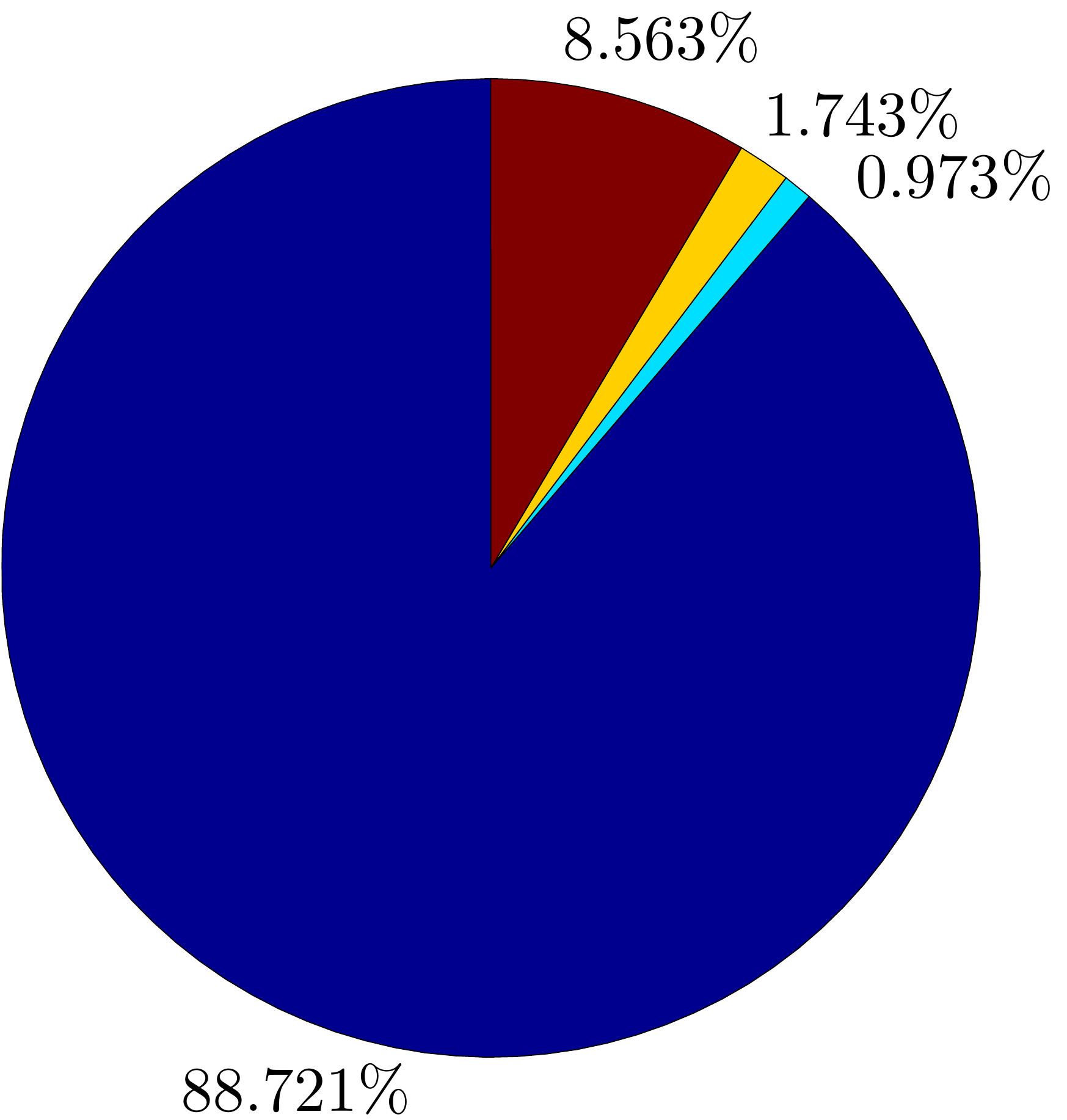}
    \caption{mult\_dcop\_03, pre-scaled.}
    \label{F:pie:noscaling:mult:equib}
\end{subfigure}
\caption{Number of possible absolute errors from dot products
in Algorithm~1 in orthogonalization kernel.
Class 1: $err < 1.0$ (blue),
Class 2: $1.0 \geq err \leq \TwoNorm{\Mat{A}}$ (light blue),
Class 3: $\TwoNorm{\Mat{A}} > err$ (yellow), and
Class 4: Non-numeric
(red).}
\label{F:pies}
\end{figure*}

We ran Algorithm~\ref{jje:alg:gmres} for 1000 total iterations, using a restart
value of 25. By instrumenting the code, we determined the numerical range of
values each vector contained, and then computed the possible absolute error that
a bit flip could introduce. We classified the absolute error according to
\S~\ref{SS:Classes}, and counted each class of errors for the duration
of the algorithm. Figure~\ref{F:pies} shows how these errors map to our
classes of errors when the matrices are scaled versus not scaled.

A large proportion of the absolute errors possible in
orthogonalization fall into Class 1 (undetectable and small). We can
explain this distribution given that the vectors $\Vect{q}_i$ are
normalized (a side effect of GMRES being derived from the Arnoldi
process).  Given normalized vectors, we know that of all the dot
products in Gram-Schmidt orthogonalization, at least one of the
vectors has data in the interval $[0,1]$. We previously established
that the interval $[0,1]$ aids in minimizing absolute error if a bit
perturbs a dot product. Now, we show how equilibrating the input
matrices can assist in forcing the non-normalized vector
($\Vect{v}_{j+1}$) as close as possible to being in the normalized
interval.

The results show the benefits of using well-scaled matrices.
Figures~\ref{F:pie:noscaling:poisson}~and~\ref{F:pie:noscaling:mult:equib}
show the Poisson problem (with no equilibration) and the
mult\_dcop\_03 matrices, which both have good scaling (see
Table~\ref{jje:table:gmres:sample_matrices:norms}).  These problems
experience a higher distribution of absolute errors less than one than
the poorly scaled CoupCons3D matrix (see
Figure~\ref{F:pie:noscaling:coup}). For the matrices that can be
equilibrated, we see that scaling the input matrices is never
detrimental, and will only improve fault tolerance, e.g., compare
CoupCons3D before scaling in Figure~\ref{F:pie:noscaling:coup} versus
after scaling in Figure~\ref{F:pie:noscaling:coup:equib}.

The pie charts are not probabilistic, that is, they do not convey the
likelihood of observing such an error. Rather, these charts
characterize the possible errors when given specific data. Consider an
arbitrary length vector $\Vect{x}$, we can determine the range of
values in the vector, e.g., $x_i \in [a,b]$, but we do not know how
many of each value, or in what order they occur. Obtaining
fine-grained statistics would involve evaluating every element of the
vector, or constructing a probabilistic model that captures the
distribution of values in each vector.

Since we consider the impact of a single bit flip, it is sufficient to follow
the methodology presented in \S~\ref{S:fault_model_eval}. That is, we may
not know the distribution and order of numbers in the vectors, but we can model
every possible error by assuming that each value in the interval \emph{could}
be used in an operation with \emph{every} value of the other interval. This
Cartesian product guarantees that we have counted all
possible errors for IEEE-754 double precision numbers in an interval, including
errors that may not occur because the vector does not contain that specific
number, or because of the ordering.

\subsubsection{Error Distribution}
Our results show that scaling tends to produce a distribution of
absolute error that is roughly $91\%$ less than or equal to one, while
$9\%$ are non-numeric.
This is expected when \emph{most} of the numbers are
near one.
Flipping the most significant exponent bit produces $11111111111$,
which will generate a non-numeric value. Similarly, the $10$
remaining exponent bits will produce error less than one --- that is,
$1/11\approx9\%$ and $10/11 \approx 91\%$. As previously discussed,
the mantissa errors are determined entirely by the exponent bits.

\subsection{Multiple Bit Flips}
While this work intentionally focuses on single bit flips, the key
finding that normalized data is \emph{better} to operate on when
performing dot products, gives some insight into how multiple bit
flips in data may behave. For example, we know that a fault in the
fractional component of a floating point number will produce an
absolute error bound above by the order of magnitude of the original
value. That is, we could flip all 52 bits of the mantissa and the
error bound from our model would still be valid (e.g.,
Eq.~\eqref{jje:eq:fault_model:dotprod:sign_bound} or
Eq.\eqref{jje:eq:magnitude_characterizes_mantissa}). In regard to
exponent flips, we have shown both analytically and experimentally
that when operating on normalized values, only 1 exponent bit per 64
bit value can introduce large error. Should the values all be
normalized, flipping $1\to0$ will minimize the value subject to
Table~\ref{jje:table:bit_flips:scalar_error}.  Experiencing more than
a single bit flip would only serve to ``shrink'' the value even
more. We intentionally do not speculate about how and why multiple bit
flips can occur, but we have shown that operating on normalized values
skews the probability of experiencing large error.

\section{Conclusion} \label{S:conclusion}

Our results indicate a clear benefit to good scaling.  We have shown
that a widely used numerical method (the Arnoldi process coupled with
Gram-Schmidt orthogonalization) inherently minimizes absolute error in
dot products. Furthermore, standard matrix equilibration algorithms
can be used to scale input matrices, which further enhance the
inherent robustness of the Arnoldi process. We demonstrated our
theoretical finding experimentally by instrumenting the GMRES
iterative solver, which is based on the Arnoldi process.

We cannot enforce that data are always normalized.  Some linear
systems may be inherently poorly scaled, or it may be impractical to
equilibrate them.  We \emph{can} advocate that scaling, while
typically used to improve numerical stability and reduce the loss of
precision, can also benefit fault resilience.  We have shown that this
result has broad applicability, because many iterative solvers are
based on orthogonal projections using normalized vectors, i.e., they
create an orthonormal basis. While this work does not propose an
end-to-end solution to soft errors, it does indicate that data scaling
can help mitigate the impact of such errors should they occur.


\bibliographystyle{IEEE}
\bibliography{paper}

\end{document}